\documentclass{article}
\usepackage{amssymb,amsmath,theorem,euscript}

\usepackage[cp1251]{inputenc}
\usepackage[english,russian]{babel}

\input{epsf}

\newcounter{sec}

\def\sm{\smallskip}

\newcounter{punct}[sec]

\def\punct{\refstepcounter{punct}{\arabic{sec}.\arabic{punct}.  }}

\def\COUNTERS{\addtocounter{sec}{1}
              \setcounter{punct}{0}
          \setcounter{equation}{0}
          \setcounter{theorem}{0}
                  }

\newtheorem{theorem}{Theorem}[sec]
\newtheorem{proposition}[theorem]{Proposition}
\newtheorem{lemma}[theorem]{Lemma}

\begin{document}

 \def\ov{\overline}
\def\wt{\widetilde}
 \newcommand{\rk}{\mathop {\mathrm {rk}}\nolimits}
\newcommand{\Aut}{\mathop {\mathrm {Aut}}\nolimits}
\newcommand{\Out}{\mathop {\mathrm {Out}}\nolimits}
 \newcommand{\tr}{\mathop {\mathrm {tr}}\nolimits}
  \newcommand{\diag}{\mathop {\mathrm {diag}}\nolimits}
  \newcommand{\supp}{\mathop {\mathrm {supp}}\nolimits}
  \newcommand{\indef}{\mathop {\mathrm {indef}}\nolimits}
  \newcommand{\dom}{\mathop {\mathrm {dom}}\nolimits}
  \newcommand{\im}{\mathop {\mathrm {im}}\nolimits}
 
\renewcommand{\Re}{\mathop {\mathrm {Re}}\nolimits}

\def\Br{\mathrm {Br}}

\def\SL{\mathrm {SL}}
\def\SU{\mathrm {SU}}
\def\GL{\mathrm {GL}}
\def\U{\mathrm U}
\def\OO{\mathrm O}
 \def\Sp{\mathrm {Sp}}
 \def\SO{\mathrm {SO}}
\def\SOS{\mathrm {SO}^*}
 \def\Diff{\mathrm{Diff}}
 \def\Vect{\mathfrak{Vect}}
\def\PGL{\mathrm {PGL}}
\def\PU{\mathrm {PU}}
\def\PSL{\mathrm {PSL}}
\def\Symp{\mathrm{Symp}}
\def\End{\mathrm{End}}
\def\Mor{\mathrm{Mor}}
\def\Aut{\mathrm{Aut}}
 \def\PB{\mathrm{PB}}
 \def\cA{\mathcal A}
\def\cB{\mathcal B}
\def\cC{\mathcal C}
\def\cD{\mathcal D}
\def\cE{\mathcal E}
\def\cF{\mathcal F}
\def\cG{\mathcal G}
\def\cH{\mathcal H}
\def\cJ{\mathcal J}
\def\cI{\mathcal I}
\def\cK{\mathcal K}
 \def\cL{\mathcal L}
\def\cM{\mathcal M}
\def\cN{\mathcal N}
 \def\cO{\mathcal O}
\def\cP{\mathcal P}
\def\cQ{\mathcal Q}
\def\cR{\mathcal R}
\def\cS{\mathcal S}
\def\cT{\mathcal T}
\def\cU{\mathcal U}
\def\cV{\mathcal V}
 \def\cW{\mathcal W}
\def\cX{\mathcal X}
 \def\cY{\mathcal Y}
 \def\cZ{\mathcal Z}
\def\0{{\ov 0}}
 \def\1{{\ov 1}}
 \def\frA{\mathfrak A}
 \def\frB{\mathfrak B}
\def\frC{\mathfrak C}
\def\frD{\mathfrak D}
\def\frE{\mathfrak E}
\def\frF{\mathfrak F}
\def\frG{\mathfrak G}
\def\frH{\mathfrak H}
\def\frI{\mathfrak I}
 \def\frJ{\mathfrak J}
 \def\frK{\mathfrak K}
 \def\frL{\mathfrak L}
\def\frM{\mathfrak M}
 \def\frN{\mathfrak N} \def\frO{\mathfrak O} \def\frP{\mathfrak P} \def\frQ{\mathfrak Q} \def\frR{\mathfrak R}
 \def\frS{\mathfrak S} \def\frT{\mathfrak T} \def\frU{\mathfrak U} \def\frV{\mathfrak V} \def\frW{\mathfrak W}
 \def\frX{\mathfrak X} \def\frY{\mathfrak Y} \def\frZ{\mathfrak Z} \def\fra{\mathfrak a} \def\frb{\mathfrak b}
 \def\frc{\mathfrak c} \def\frd{\mathfrak d} \def\fre{\mathfrak e} \def\frf{\mathfrak f} \def\frg{\mathfrak g}
 \def\frh{\mathfrak h} \def\fri{\mathfrak i} \def\frj{\mathfrak j} \def\frk{\mathfrak k} \def\frl{\mathfrak l}
 \def\frm{\mathfrak m} \def\frn{\mathfrak n} \def\fro{\mathfrak o} \def\frp{\mathfrak p} \def\frq{\mathfrak q}
 \def\frr{\mathfrak r} \def\frs{\mathfrak s} \def\frt{\mathfrak t} \def\fru{\mathfrak u} \def\frv{\mathfrak v}
 \def\frw{\mathfrak w} \def\frx{\mathfrak x} \def\fry{\mathfrak y} \def\frz{\mathfrak z} \def\frsp{\mathfrak{sp}}
 \def\bfa{\mathbf a} \def\bfb{\mathbf b} \def\bfc{\mathbf c} \def\bfd{\mathbf d} \def\bfe{\mathbf e} \def\bff{\mathbf f}
 \def\bfg{\mathbf g} \def\bfh{\mathbf h} \def\bfi{\mathbf i} \def\bfj{\mathbf j} \def\bfk{\mathbf k} \def\bfl{\mathbf l}
 \def\bfm{\mathbf m} \def\bfn{\mathbf n} \def\bfo{\mathbf o} \def\bfp{\mathbf p} \def\bfq{\mathbf q} \def\bfr{\mathbf r}
 \def\bfs{\mathbf s} \def\bft{\mathbf t} \def\bfu{\mathbf u} \def\bfv{\mathbf v} \def\bfw{\mathbf w} \def\bfx{\mathbf x}
 \def\bfy{\mathbf y} \def\bfz{\mathbf z} \def\bfA{\mathbf A} \def\bfB{\mathbf B} \def\bfC{\mathbf C} \def\bfD{\mathbf D}
 \def\bfE{\mathbf E} \def\bfF{\mathbf F} \def\bfG{\mathbf G} \def\bfH{\mathbf H} \def\bfI{\mathbf I} \def\bfJ{\mathbf J}
 \def\bfK{\mathbf K} \def\bfL{\mathbf L} \def\bfM{\mathbf M} \def\bfN{\mathbf N} \def\bfO{\mathbf O} \def\bfP{\mathbf P}
 \def\bfQ{\mathbf Q} \def\bfR{\mathbf R} \def\bfS{\mathbf S} \def\bfT{\mathbf T} \def\bfU{\mathbf U} \def\bfV{\mathbf V}
 \def\bfW{\mathbf W} \def\bfX{\mathbf X} \def\bfY{\mathbf Y} \def\bfZ{\mathbf Z} \def\bfw{\mathbf w}
 \def\R {{\mathbb R }} \def\C {{\mathbb C }} \def\Z{{\mathbb Z}} \def\H{{\mathbb H}} \def\K{{\mathbb K}}
 \def\N{{\mathbb N}} \def\Q{{\mathbb Q}} \def\A{{\mathbb A}} \def\T{\mathbb T} \def\P{\mathbb P} \def\G{\mathbb G}
 \def\bbA{\mathbb A} \def\bbB{\mathbb B} \def\bbD{\mathbb D} \def\bbE{\mathbb E} \def\bbF{\mathbb F} \def\bbG{\mathbb G}
 \def\bbI{\mathbb I} \def\bbJ{\mathbb J} \def\bbK{\mathbb K} \def\bbL{\mathbb L} \def\bbM{\mathbb M} \def\bbN{\mathbb N} \def\bbO{\mathbb O}
 \def\bbP{\mathbb P} \def\bbQ{\mathbb Q} \def\bbS{\mathbb S} \def\bbT{\mathbb T} \def\bbU{\mathbb U} \def\bbV{\mathbb V}
 \def\bbW{\mathbb W} \def\bbX{\mathbb X} \def\bbY{\mathbb Y} \def\kappa{\varkappa} \def\epsilon{\varepsilon}
 \def\phi{\varphi} \def\le{\leqslant} \def\ge{\geqslant}

\def\UU{\bbU}
\def\Mat{\mathrm{Mat}}
\def\tto{\rightrightarrows}

\def\Gr{\mathrm{Gr}}

\def\graph{\mathrm{graph}}

\def\O{\mathbb{O}}

\def\la{\langle}
\def\ra{\rangle}

\def\B{\mathrm B}
\def\Int{\mathrm{Int}}
\def\LGr{\mathrm{LGr}}


 \newcommand{\ess}{\mathop {\mathrm {ess\,sup}}\nolimits}

\def\sfX{\mathsf X}
\def\sfY{\mathsf Y}
\def\sfZ{\mathsf Z}
\def\sfV{\mathsf V}
\def\sfU{\mathsf U}
\def\sfW{\mathsf W}

\def\I{\mathbb I}
\def\M{\mathbb M}
\def\T{\mathbb T}

\def\Ams{\mathrm{Ams}}
\def\Gms{\mathrm{Gms}} 
\def\Mar{\mathrm{Mar}}
\def\Pol{\mathrm{Pol}}
\def\Naz{\mathrm{Naz}}
\def\naz{\mathrm{naz}}
\def\bNaz{\mathbf{Naz}}
\def\AMod{\mathrm{AMod}}
\def\ALat{\mathrm{ALat}}

\def\Ver{\mathrm{Vert}}
\def\Bd{\mathrm{Bd}}
\def\We{\mathrm{We}}
\def\Heis{\mathrm{Heis}}

\def\bbot{{\bot\!\!\!\bot}}
\def\tri{\triangledown}
\def\zigzag{\rightsquigarrow}

\def\lk{<<}
\def\rk{>>\,}
\def\rkk{>>}

\begin{flushright}
 UDC 517.518.112, 512.583,  517.983.23
\end{flushright}

\begin{center}
\Large\bf

On the boundary of the group of transformations leaving a measure quasi-invariant 

\bigskip

\large\sc
Yu.A.Neretin%
\footnote{Supported by the grant FWF,  Project 22122, and RosAtom, Contract H.4e.45.90.11.1059.
}
\end{center}


{\small Let $A$ be a Lebesgue measure space. We interpret measures on 
$A\times A\times \R^\times$  as 'maps'   $A$ to $A$, which 'spread' 
  $A$ along itself;
  their Radon-Nikodym derivatives also are spread.
We discuss basic properties of the semigroup of such maps
and action of this semigroup in spaces  $L^p(A)$.}

\section{Purposes of the work}

\COUNTERS

{\bf\punct Groups $\Ams(A)$, $\Gms(A)$ and their boundaries.} 
Denote by  
$\R^\times$ the multiplicative group of positive
real numbers.
Let   $A$ be a space with continuous probabilistic measure  $\alpha$.
Denote by $\Ams(A)$  the group of measurable transformations of the space
  $A$ preserving $\alpha$,
 by $\Gms(A)$ we denote the group of transformations leaving the measure
$\alpha$ quasiinvariant.

The group $\Ams(A)$ has a well-known completion  $\ov{\Ams(A)}$ 
(below we denote it by  $\Mar(A,A)$), points of the completion
are measures on $A\times A$ whose projections to both factors coincide with
  $\alpha$.
Elements of $\ov{\Ams(A)}$  can be regarded as 'maps'
 $A\to A$ spreading points along the set
  $A$. There is a well-defined composition of such maps. 

Such objects are widely used in probability
(this is simply a reformulation of notion of
'Markov operators') and in ergodic theory
 (see, e.g., \cite{Kre}, \cite{Ver}, \cite{Rud}, \cite{Gla}),
 they appear in mathematical hydrodynamics
(see, e.g.,  \cite{Bre}).

The group $\Gms(A)$ also has a natural completion   $\ov{\Gms(A)} $ 
(below we denote it by  $\Pol(A,A)$), 
whose points are measures
on  $A\times A\times \R^\times$, 
such measures can be regarded as spreading maps
with spread Radon--Nikodym derivative;
we call such 'maps' by 'polymorphisms'%
\footnote{May by, it is better to say   '$\R^\times$-polymorphisms'.
A.M.Vershik \cite{Ver} introduced the term 'polymorphism'
for elements of   $\Mar$,
more common is the term 'bistochastic kernels'.
 In \cite{Ner-bist}, \cite{Ner-book}, \cite{Ner-discr} 
there were considered measures on 
   $A\times A\times G$, where  $G$
   is an arbitrary group, they were called $G$-polymorphisms.}%
$^,$%
\footnote{This objects differ from 'substochastic kernels' \cite{Kre}.}.
The semigroup  $\ov{\Gms(A)} $ was introduced in
 \cite{Ner-bist},
 and an initial
motivation was the following theorem:

{\it Any unitary representation of the group $\Gms(A)$ admits a unique continuous
extension to the semigroup  $\ov{\Gms(A)}$.}

\sm

{\bf\punct Spaces $\Pol(A,B)$ and  multiplication.}
Denote by
  $t$ the coordinate on  $\R^\times$. Denote by  $\cM$
  the semigroup of positive finite measures on the group
  $\R^\times$.

Consider two Lebesgue spaces with measure
 $(A,\alpha)$, $(B,\beta)$.
 We say that a measure on
 $A\times B$ is a {\it polymorphism} 
if

\sm

1)  The pushforward of $\frP$  under the projection
to the first factor  $A$ coincides with    $\alpha$.

\sm

2) The pushforward of  $t\times \frP$
under the projection to the second factor
 $B$ coincides with $\beta$.

\sm

Denote the set of all such measures by
 $\Pol(A,B)$.

 We embed the group $\Gms(A)$
to $\Pol(A,A)$ by the following rule.
Let  $g\in \Gms(A)$. Consider the map 
$A\to A\times A\times \R^\times$
given by the formula
$$
a\mapsto \bigl(a,g(a),g'(a)\bigr),
$$
where $g'(a)$ denote the Radon--Nikodym derivative.
The image of the measure  $\alpha$ under this map
is a measure on
 $A\times A\times \R^\times$, which satisfies the properties
 1) and 2).

It turns out that there exists a natural multiplication
 $$
 \Pol(A,B)\times \Pol(B,C)\to \Pol (A,C).
 $$
 It looks transparently in the following special case.
Consider a map  $p$, which  for almost  each
 $(a,b)\in A\times B$ assigns a  measure 
$p_{a,b}$ on $\R^\times$. 
Then we can define a measure  $\frP$ on $A\times B\times \R^\times$
from the following condition.
Let $M\subset A$, $N\subset B$, $K\subset \R^\times$ be measurable sets.
Then
$$
\frP(M\times N\times \R^\times)=\int_{M\times N} p_{a,b}(K)\,dp(a,b)
,
$$
where $p_{a,b}(K)$ is a measure $p_{a,b}$ of the set  $K$.
The measure must satisfy the properties
 1), 2),  this implies sufficiently evident conditions
 for  $p$
(see below  (\ref{eq:abs1}), (\ref{eq:abs2})).
Take two maps
 $p:A\times B\to \cM$, $q: B\times C\to\cM$
and the corresponding polymorphisms  $\frP$, $\frQ$.
Consider the  map
 $r:A\times C\to\cM$ given by the formula 
$$
r=\int_B p_{a,b}*q_{b,c}\,d\beta(b),
$$
where $*$ denotes the convolution of measures
on the multiplicative group  $\R^\times$. 
The corresponding polymorphism
 $\frR$ is a product  $\frP\circ \frQ$.

It turns out  that this product can be extended  by
a separate continuity
 (see Subsection \ref{ss:pol-convergence}) to an associative
 operation on arbitrary polymorphisms. One of the purposes
  of the work is to give
various operational definitions for this product%
 \footnote{Existence of this product is not self-obvious.
 A detailed written proof is tiresome,
 the argumentation of
 \cite{Ner-bist}  uses a dual language 
 (Theorem \ref{th:pol-mellin-product} below can be regarded as a definition
 of the product), but this way also is long.
 According the kernel theorem 
 (see, e.g., \cite{Her}),
 any operator in $L^2(\R)$ is an integral operator
 in the sense of L.~Schwartz. However, a calculation
 of the kernel of a product by the usual formula 
  \begin{equation}M(x,z)=\int K(x,y)L(y,z)\,dy,\label{eq:convolution-kernels}\end{equation}
  generally fails (the integral can diverge, even integrand can be not well defined).
In our case, a value of a singular expression of the form
 (\ref{eq:convolution-kernels}) can be defined.}.
 
 As a result, we get a category, whose objects 
 are Lebesgue spaces,  morphisms are 
 are polymorphisms.

\sm

{\bf \punct  Action in the spaces  $L^p$.} 
 Let $u=v+iw$ ranges in the strip 
$0\le v\le 1$ in $\C$.
For any 
 $g\in \Gms(A)$ we assign a family  $u\mapsto T_{u}(g)$  of linear operators
 in the space of measurable functions on
 $A$ by
$$
T_{v+iw} f(g)=f(g(a))g'(a)^{v+iw}
.
$$
Evidently,
an operator $T_{v+iw}$ is an isometry of the space
$L^{1/v}$. Thus we get a family of representations
of the group   $\Gms(A)$ depending holomorphically on the parameter
$u$.

It turns out that the representations  $T_{u}$
can be extended to the category of polymorphisms,
namely for any polymorphism $\frP\in \Pol(A,B)$
there is family of linear operators
$$
T_u(\frP):L^{1/v}(B)\to L^{1/v}(A)
$$
 such that
$$
T_u(\frP\circ \frQ)=T_u(\frQ)T_u(\frP)
$$
и
$$
\|T_{v+iw}(\frP)\|_{L^{1/v}}\le 1
.$$

For polymorphisms defined by a function
  $p$ as above the operator $T_u(\frP)$ equals
$$
T_u(\frP) f(a)=\int_{B} \int_{\R^\times} f(b)t^u dp_{a,b}(t)\, d\beta(b).
$$

\sm


{\bf\punct Olshanski's problem on weak closure.} Let $\rho$ 
be a unitary representation of a group
 $G$ in a Hilbert space   $H$. Consider the set   $\rho(G)$
 of all operators
$\rho(g)$, where $g$ ranges in  $G$. Consider its closure 
 $\Gamma=\Gamma_\rho=\ov{\rho(G)}$ with respect to the
 weak operator topology.
It easy to show that   $\Gamma$
is a compact semigroup. For Lie groups this object
is not interesting 
(usually we get a one-point compactification
 $G$, see \cite{HM}). For infinite-dimensional groups the picture
 changes.
 The following 'experimental facts' hold
(see \cite{Olsh-topics}, \cite{Ner-book}). 

\sm

--- The semigroup $\Gamma_\rho$ is essentially larger than    $G$.

\sm

--- $\Gamma=\Gamma_\rho$ admits a universalization  
  (a {\it mantle} of group $G$) with respect to  $\rho$.

\sm

--- $\Gamma$ admits an explicit description.

\sm

--- it turns out that  $\Gamma$ is an effective tool
for investigation of representations of the group
 $G$.

\sm


{\bf \punct Action of the mantle on measure spaces.}
Consider an action of infinite-dimensional group 
 $G$ by transformations leaving the measure quasiinvariant 
  (many such actions are known, see survey
     \cite{Ner-frac}
     and relatively recent constructions 
 \cite{Ner-hua}, \cite{KOV}, \cite{BO}, \cite{Ner-p}, \cite{Ner-finite}).  In \cite{Ner-frac}
 there was proposed arguments (partially formal, partially heuristic),
 which show that the  mantle 
 $\Gamma$ acts on   $A$ by polymorphisms.

In \cite{Ner-match}, \cite{Ner-coll} such actions were
described in two simplest cases:
for groups of natural symmetries of Gaussian
and Poisson measures. By  author's opinion, formulas
looks unusual.  Therefore there arises a problem about description
of such actions in more complicated cases.
This problem can be reformulated in spirit of Olshanski:
to describe the closure of
  $G$ in $\ov{\Gms(A)}$.

\sm


{\bf\punct  Purposes of the paper.} Several statements about polymorphisms were
formulated in
 \cite{Ner-bist}, \cite{Ner-match} without proofs. 
 The present paper is a step backward, here we present proofs,
 this also underpin the papers \cite{Ner-match},
\cite{Ner-coll} and the problem formulated above.
In the paper we give several equivalent definitions of
the  product of polymorphisms,
prove their self-consistency,
and describe the correspondence between polymorphisms and holomorphic 
operator-valued functions in the strip.
This provides a dual language for work with polymorphisms
(see \cite{Ner-match}, \cite{Ner-coll}),
in particular, we get a non-direct but convenient definitions 
of the  product of polymorphisms.


\sm


{\bf\punct  The structure of the paper.}
Sections  2 and 3 contain preliminaries
on Lebesgue spaces and Markov operators.
In Section
 4 we describe some properties of the semi-ring
of positive measures on  
 $\R^\times$.
 Polymorphisms are defined
in \S 5. In \S 6
we discuss action of polymorphisms on spaces
 $L^p$.

\sm



\section{Preliminaries. Lebesgue spaces}

\COUNTERS

This section contains several standard definitions and notations, which are used below.
A fundamental work on Lebesgue spaces is the paper by Rokhlin
 \cite{Roh}.
Its exposition is contained in \cite{MI}.

\sm


{\bf\punct Lebesgue spaces.%
\label{ss:lebesgue-spaces}}
 A {\it Lebesgue space}%
 \footnote{Some authors use a more precise term 'Lebesgue--Rokhlin space'.
 The term 'Lebesgue space' is ambiguous but more generally accepted.}
  $(A,\alpha)$
is a space with a positive finite measure
equivalent to a disjoint union of a finite segment
  $[p,q]\subset \R$ with Lebesgue measure and a finite or countable
  collection of points (atoms) having non-zero measure.
We assume  $\alpha(A)>0$.

A measure is called

\sm

--- {\it probabilistic} if $\alpha(A)=1$;

\sm

--- {\it continuous} if the set of atoms is empty;

\sm

--- {\it discrete} if  $A$ is a union of atoms.

\sm

It is known that almost all spaces with a {\it finite}
Lebesgue measure that arise in analysis are Lebesgue.

\sm

Such a space (a union of a segment and a collection of atoms)
has a natural Borel structure. Below the term {\it measurable set}
denotes a measurability with respect to a Borel structure,
the term 
 {\it measure} means a measure defined on a Borel  $\sigma$-algebra.

We denote by
 $\alpha(M)$ the  measure of a measurable set 
 $M\subset A$.
By $\int f(a)\,d\alpha(a)$ we denote the integral with respect
to a measure 
$\alpha$.


\sm

{\bf \punct Spaces $L^p$.%
\label{ss:Lp}}
Let $1\le p<\infty$. Consider the space   $L^p(A)$
consisting of measurable functions $f$
(defined upto a.s.)
 satisfying the condition 
$$
\|f\|_p:=\biggl(\int_A|f(a)|^p\,d\alpha(a)\biggr)^{1/p}<\infty
.$$
We get a separable Banach space with norm
 $ \|f\|_p$. 
If $p>r$, then  $L^p(A)\subset L^r(A)$.

For $p=\infty$ we define the norm as%
\footnote{
Recall that the {\it essential supremum} of a set
 $X\subset \R$ is the infimum of all  
$x$ such that the measure of  $X\cap [x,\infty)$ is 0.}
$$
\|f\|_\infty:=\ess\limits_{a\in A}|f(a)|
.
$$
In this case we get a nonseparable Banach space.
To escape the non-separability%
\footnote{We wish to use a duality, but the space
dual to
 $L^\infty$ is a pathological object, see, e.g., 
 \cite{KA}, \S.IV.2},  
we introduce the space  $L^{\infty-}(A)$.
As above it consists of bounded measurable function
 {\it but we change a definition of convergence in 
 $L^\infty$}. We assume that a sequence   $f_j\in L^\infty$
converges to  $f$ if the sequence 
 $\|f_j\|_\infty$  is bounded and  for any $\epsilon>0$
 the measure of the set
 $\bigl\{a\in A:\,|f_j(a)-f(a)|>\epsilon\bigr\}$ tends to $0$ 
as $j$ tends to  $\infty$. 
We say that  a linear functional  $\ell$ is continuous 
$L^{\infty_-}(A)$ if it is bounded in the sense
of $L^\infty$ and the convergence
 $f_j\to f$ implies the convergence  $\ell(f_j)\to\ell(f)$. 

\begin{lemma}
\label{l:infty=p}
On bounded with respect to $L^\infty(A)$-norm subsets,  $L^{\infty_-}$-convergence
is equivalent to $L^p$-convergence for any   $p<\infty$,
and also is equivalent to the convergence in measure.
\end{lemma}

Proof is obvious.

\sm

Let $\frac 1p+\frac 1q=1$, $p\ne \infty$.
Any continuous  linear functional
on  $L^p(A)$ has the form 
$$
\gamma(f)=\int_A f(a) g(a)\,d\alpha(a), \qquad\text{where $g\in L^q(A,\alpha)$}
.$$
Moreover, $\|\gamma\|=\|g\|_q$.

\begin{lemma}
Any continuous linear functional on $L^{\infty_-}(A)$
has the same form with $g\in L^1(A)$.
\end{lemma}

{\sc Proof.} 
We evaluate  $\gamma$ on a characteristic function
of a measurable set and get a countably additive
charge on 
 $A$. For sets of zero measure this charge is   $0$.
By the Radon--Nikodym theorem (see \cite{KF}) 
this charge is determined by a measurable integrable function. 
 \hfill $\square$
 
 \sm

{\sc Remark.} The convergence in $L^{\infty_-}(A)$
corresponds to a locally convex topology determined by
the following family of seminorms. First,
for any
 $h\in L^1(A)$ we define a seminorm 
$$
[f]_h:=\Bigl| \int_A f(a) h(a)\,d\alpha(a)\Bigr|
.
$$
To the family $[f]_h$ (which determines an $L^1$-weak
topology on  $L^\infty(A)$) we add
the norm $\|f\|_{L^1}$. We will not use this.
\hfill $\square$

\sm


{\bf \punct Image of  measure.%
\label{ss:pushforward}}
Let $(A,\alpha)$  be a Lebesgue space, $B$ a space with the standard
Borel structure.
Consider a measurable map
 $\pi:A\to B$.
The measure  $\beta$ on $B$ is defined by the condition:
$\beta(N)=\alpha\bigl(\pi^{-1}(N)\bigr)$.
The space 
 $(B,\beta)$ obtained in this way is Lebesgue.

\sm


{\bf \punct Conditional measures.%
\label{ss:conditional-measures}}
A countable (or finite) partition 
 $\sfX$ of a Lebesgue space 
 $(A,\alpha)$  is a representation of  
$A$ as a disjoint union of measurable sets, $\sfX:A=\cup X_j$.
The quotient  $A/\sfX$ is a discrete space consisting of points
  $a_j$ with measures
$\alpha(X_j)$.

A continual {\it partition} $\sfX:A=\cup_{r\in R} X_{r}$, 
where $r$ range a continual space  $R$ and $X_r$ are mutually disjoint,
is called   {\it measurable}%
\footnote{See \cite{Roh}, \cite{MI}. A partition of $\R$ into equivalence classes
$x\sim y$ if $x-y\in \Q$ is an example of a non-measurable partition.}
if there is a countable family of measurable subsets
  $U_j\subset A$
  such that

\sm

--- any $U_j$ is a union  $\cup_{r\in P}  X_r$, where  $P\subset R$ is a subset;

\sm

--- the family $U_j$ separates
 $X_r$,  i.e., for any $X_r\ne X_q$
there is   $U_i$ such that   $X_r\subset U_i$, $X_q\not\subset U_i$.

\sm

We define a structure of a measure spaces on the quotient
 $A/\sfX\simeq R$: a subset
$P\subset R$ is measurable, iff   $\cup_{r\in P} X_r$ is measurable,
and measure
of
 $P$ is $\rho(P):=\alpha\bigl(\cup_{r\in P} X_r\bigr)$.

The space
 $A/\sfX$ obtained in this way is Lebesgue, the map  $A\to A/\sfX$ 
 is measurable.

Conversely, for any measurable map of Lebesgue spaces
 $g:A\to B$ the partition
$A=\cup_{b\in B} g^{-1}(b)$ is measurable.

\sm

Recall the Rokhlin theorem.
{\it For any measurable partition   $\sfX:A=\cup_{r\in R} X_r$
there exists a family of probability measures
$\xi_r$ defined for almost all  
 {\rm(}with respect to the measure on $A/\sfX${\rm)}  sets  $X_r$
 such that for any measurable subset
  $M\subset A$ and for almost all   $r\in R$ the subsets 
$M\cap X_r\subset X_r$ are measurable in  $X_r$ and
$$
\alpha(M)=\int_{A/\sfX} \xi_r(M\cap X_r) d\rho(r)
.$$
Almost all spaces
 $X_r$ are Lebesgue.
For any integrable function on $A$ the following identity holds} 
$$
\int_A f(a)\,d\alpha(a)= \int_{A/\sfX} \int_{a\in X_r} f(a)\, d\xi_r(a)\, d\rho(r)
.$$
The measures
 $\xi_r$  are called  {\it conditional measures}.


\sm

{\bf\punct Conditional expectations.%
\label{ss:conditional-expectation}}
Let $R=A/\sfX$  be the quotient space, $\pi:A\to R$  the projection map,
 $\xi_r$ the conditional measures.
The operator {\it of conditional expectation} 
$$J[A;\sfX]:L^1(A)\to L^1(R)$$
is defined by
$$
J[A;\sfX]f(r)=\int_{X_r} f(a)\,d\xi_r(a).
$$

On the other hand there is an isometric embedding
$$K[A;\sfX]:L^1(R)\to L^1(A),
$$
defined by
$$
K[A;\sfX]h(a)=h(\pi(a))
.
$$

We also define the operator 
 {\it of conditional average}
$$
I[A;\sfX]=K[A;\sfX]\,J[A;\sfX]: L^1(A)\to L^1(A)
.
$$
It can be represented as
$$
I[A;\sfX]f(a)=\int_{X_p\ni a} f(c)\,d \xi_p(c)
.
$$
These operators satisfy the following identities
$$
I^2=I,\qquad IK=K, \qquad JI=J,\qquad JK=1
.
$$

\sm


{\bf \punct The group $\Ams(A)$.%
\label{ss:ams}}
Let $(A,\alpha)$ be a  Lebesgue  space with a continuous measure.
By $\Ams(A)$ we denote the group of measure preserving bijections (a.s.)
 $A\to A$.
 Two elements
   $g_1$, $g_2$ of the group  $\Ams(A)$ are considered as coinciding
   if 
$g_1(a)=g_2(a)$ a.s.

The group $\Ams(A)$ acts in the space   $L^p(A,\alpha)$
by isometric operators
$$
T(g)f(a)=f\bigl(g(a)\bigr).
$$
This group is a separable topological group.
The convergence is defined by the condition:
 $g_j\to g$  if for all measurable subsets   $M$, $N\subset A$ 
 we have
 $$
\lim_{j\to\infty} \alpha\bigl(g_j(M)\cap N\bigr)=\alpha\bigl(g(M)\cap N\bigr)
 .
 $$


\sm

{\bf \punct The group $\Gms(A)$.%
\label{ss:gms}}
Recall that a measure
 $\alpha$ is {\it quasiinvariant}  with respect to bijective a.s.s
 map 
$A\to A$ if for any subset  $M\subset A$ of zero measure,
sets  $g(A)$ and $g^{-1}(A)$ have zero measure.

An equivalent condition:
there is a function
 $g'(a)$, which is called   {\it Radon--Nikodym derivative} (see, \cite{KF})
 such that for any measurable set  $M\subset A$
the following equality holds 
$$
\mu(gM)=\int_M g'(a)\,d\alpha(a),
$$
and $g'(a)\ne 0$ a.s. on $A$.

The Radon-Nikodym derivative satisfies the usual chain rule:
$$
(g\circ h)'(a)=g'\bigl(h(a)\bigr)\,h'(a)
.
$$

By $\Gms(A)$ we denote the group of  bijections a.s.
  $A\to A$ leaving the measure $\alpha$ quasiinvariant.

Fix $p$. For each  $s\in\R$ we define an  action 
of the group  $\Gms(A)$ in  $L^p(A,\alpha)$ 
by isometric operators by the formula
\begin{equation}
T_{1/p+is}(g) f(a)=f\bigl(g(a)\bigr) g'(a)^{1/p+is}
.
\label{eq:t1p}
\end{equation}

According the chain rule
these operators satisfy
$$
T_{1/p+is}(g_1) T_{1/p+is}(g_2)=T_{1/p+is}(g_1\circ g_2)
.$$

\section{Markov category}

\COUNTERS

Bistochastic kernels and Markov operators
discussed below are standard objects, see, e.g.,
 \cite{Ver}, \cite{Kre}, \cite{Ner-book}, \cite{Gla}. 
 For a coherence of the text sometimes we present proofs or sketches
 of proofs.


\sm

{\bf\punct Markov category.%
\label{ss:markov-category}} Objects of the category  $\Mar$
are Lebesgue spaces with probabilistic measure.
A morphism $\frp:(A,\alpha)\to (B,\beta)$ (a {\it bistochastic kernel})
is a measure
$\frp$ on $A\times B$ such that

\sm

--- the image of  $\frp$ under the projection   $A\times B\to A$ 
coincides with   $\alpha$;

\sm

--- the image of $\frp$ under the projection  $A\times B\to B$
coincides with  $\beta$.

\sm

We denote the set of all morphisms 
 $\frp:(A,\alpha)\to (B,\beta)$ by $\Mar(A,B)$.
 A general rule for multiplication of morphisms is 
 little below
(Subsection \ref{ss:Mar-general}).
Before this we consider a transparent special case.

\sm


{\bf\punct A special case: spaces with discrete measure.%
\label{ss:mar-discrete}}
Let spaces 
 $A$, $B$ be countable. Let  $a_i$ (resp. $b_j$) be their points. 
 Denote by 
$\alpha_i$ (resp. $\beta_j$) their measures.
We can regard morphisms
 $\frp\in\Mar(A,B)$ as matrices   $\frP=\frp_{ij}$ such that
$$
\frp_{ij}\ge 0, \qquad \sum_i \frp_{ij}=\beta_j, \qquad \sum_j \frp_{ij}=\alpha_i
.$$
For $\frp\in\Mar(A,B)$, $\frq\in\Mar(B,C)$ the product is given
by the formula 
$$
\frr_{ik}=\sum_j \frac{\frp_{ij}\frq_{jk}}{\beta_j}
$$
or
\begin{equation}
\frR= \frQ \Delta_\beta^{-1}\frP
,
\label{eq:psi-delta-phi}
\end{equation}
where $\Delta_\beta$ is a diagonal matrix
with elements   $\beta_j$ on the diagonal.

\sm


{\bf\punct   A special case: absolutely continuous kernels.%
\label{ss:markov-absolute}}
Let $p:A\times B\to \R$ be a nonnegative integrable function
satisfying the conditions 
$$
\int_B p(a,b)\,d\beta(b)=1
\qquad\int_A p(a,b)\,d\alpha(a)=1\qquad\text{a.s.}
$$
Then we define a bistochastic kernel
 $\frp$ on $A\times B$
from the condition
$$
\frp(M\times N)=\int_M\int_N p(a,b) \,d\beta(b) \,d\alpha(a)
$$
Let $p:A\times B\to \R$, $q:B\times C\to \R$ be such functions.
The product of the corresponding bistochastic kernels  corresponds 
to the function
\begin{equation}
r(a,c):=\int_B p(a,b)\,q(b,c)\,d\beta(b)
\label{eq:for-fubini}
\end{equation}
(this is the usual formula for product of integral operators).

\begin{lemma}
For almost all
 $c$ and almost all  $a$ the integral converges. 
\end{lemma}

{\sc Proof.} For almost all $c$, 
$$
\int_B\int_A p(a,b)\,q(b,c)\,d\alpha(a)\,d\beta(b)=\int_B q(b,c)\,d\beta(b)=1
.
$$
Applying the Fubini Theorem
(see \cite{KF})
we observe that integral
 (\ref{eq:for-fubini})
converges for almost all 
 $a$.
\hfill$\square$


\sm

{\bf\punct Definition of the product in the general case.%
\label{ss:Mar-general}}
Consider morphisms $\frp:(A,\alpha)\to(B,\beta)$, $\frq:(B,\beta)\to (C,\gamma)$.
We wish to define their product
 $\frr=\frq\circ\frp:(A,\alpha)\to (C,\beta)$.
 Consider  $M\subset A$, $K\subset C$.
Restrict $\frp$ to $M\times B$ and consider the image   $\frp_{M,b}$ of the restriction
under the projection $M\times B\to B$.
Since $\frp_M(b)$ is dominated by  $\beta(b)$
we get
\begin{equation}
\frp_M(b)= u_M(b)\,d\beta(b),
\label{eq:u}
\end{equation}
where $u_M(b)$ is a positive function  $\le 1$.
In a similar way, consider the restriction
of  $\frq$ to
$B\times K$ and represent the image $\frq_{K}(b)$ of the restriction under the map
 $B\times K\to B$ as
$$
\frq_K(b)=
v_K(b)\,d\beta(b)
.$$
Again, $0\le v_N(b)\le 1$.
We assume
\begin{equation}
\frr(M\times K):=\int_B u_M(b)\,v_K(b)\,d\beta(b).
\label{eq:uv}
\end{equation}
 
\begin{proposition}
\label{pr:mar-associativity}
The multiplication
  $\Mar(A,B)\times\Mar(B,C)\to\Mar(A,C)$ 
defined in this way is associative.
\end{proposition}


\begin{lemma}
\label{l:sovpad}
For absolutely continuous kernels the multiplication defined
in this way coincides with the  multiplication defined in Subsection
{\rm\ref{ss:markov-absolute}}.
\end{lemma}

{\sc Proof.}  Evaluate the measure of a set  $M\times K\subset A\times B$. In 
notation of
(\ref{eq:for-fubini}),
\begin{multline*}
\int_{M\times K} r(a,c)\,d\alpha(a)\,d\gamma(c)=
\int_{M\times K} \Bigl(
\int_B p(a,b)\,q(b,c)\,d\beta(b)\Bigr) \,d\alpha(a)\,d\gamma(c)
=\\=
\int_B
\Bigl(\int_M p(a,b)\,d\alpha(a) \Bigr)\cdot\Bigl( \int_K q(b,c) d\gamma(c)\Bigr)\,\,d\beta(b)
\end{multline*}
The right-hand side coincides with
 (\ref{eq:uv}).
\hfill $\square$

\sm


{\bf\punct Involution.}
The identical map
 $A\times B\to B\times A$ induces the map   $\Mar(A,B)\to\Mar(B,A)$. 
We denote it by $\frp\mapsto \frp^\bigstar$. Obviously,
$$
(\frq\circ\frp)^\bigstar=\frp^\bigstar\circ\frq^\bigstar
.
$$


\sm


{\bf\punct Automorphisms.%
\label{ss:automorphism}} Let $(A,\alpha)$ be a space with continuous 
measure.
Let $g\in\Ams(A)$. Consider a map 
$\iota_g:A\to A\times A$ defined by the formula 
 $\iota(a)=\bigl(a,g(a)\bigr)$.
Denote by
 $\xi[g]$ the image of the measure   $\alpha$ under this
 map.
 Evidently, $\xi[g]\in\Mar(A,A)$.
 It is easy to see that
$$
\xi[g_1g_2]=\xi[g_1]\xi[g_2].
$$


{\bf\punct Convergence.%
\label{ss:mar-convergence}} A sequence $\frp_j\in \Mar(A,B)$
converges to
$\frp\in \Mar(A,B)$ if for any subsets
$M\subset A$, $N\subset B$ we have
$$
\lim_{j\to\infty} \frp_j(M\times N)=  \frp(M\times N)
.$$

\begin{proposition}
\label{pr:mar-product}
{\rm a)} The space $\Mar(A,B)$ is metrizable and compact.

\sm

{\rm b)} Absolutely continuous measures are dense in $\Mar(A,B)$.

\sm

{\rm c)} The product $\Mar(A,B)\times\Mar(B,C)\to\Mar(A,C)$
is separately continuous.
\end{proposition}

{\sc Proof.}
a) We consider a countable family of sets
  $M_i\subset A$ such that for any
 $M\subset A$ and any   $\epsilon>0$ there exists
$M_i$ such that
$$
\alpha\bigl((M\setminus M_i)\cup (M_i\setminus M)\bigr)<\epsilon.
$$
In a similar way, we choose a family
 $N_j\subset N$.  A metric is defined by 
$$
\rho(\frp,\frq)=\sum_{i,j} 3^{-i-j} \bigl|\frp(M_i\times N_j)-\frq(M_i\times N_j)\bigr|
.
$$
Compactness can be proved by the usual diagonal procedure
 (see \cite{RS1}, Theorem   I.24).

\sm

c) A convergence $\frp_j$ to $\frp$ implies the weak convergence 
 (see \cite{RS1}, \cite{KF}) of functions
$u_j\to u$ (see (\ref{eq:u})) in the sense of $L^2(B)$.
Formula (\ref{eq:uv}) has a form of an inner product,
and the inner product is separately continuous with respect to the weak
convergence.
\hfill $\square$ 

\sm

{\sc Proof of Proposition  \ref{pr:mar-associativity} 
(associativity of multiplication).}
For absolutely continuous kernels the associativity is obvious,
it remains to refer to separate continuity of the multiplication.

\hfill $\square$

\begin{proposition}
\label{pr:ams-density}
Let a measure
 $\alpha$ be continuous. Then the group   $\Ams(A)$ is dense in $\Mar(A,A)$.
\end{proposition}

See \cite{Ner-book}, Theorem 4.4.1.



\sm

{\bf\punct Another language and an equivalent definition of
the product.%
\label{ss:mar-def-2}}
For $\frp\in\Mar(A,B)$ consider the projection    $A\times B\to B$. 
We have conditional measures
 $\frp_a(b)$ on almost all fibers, they satisfy 
 the equation 
$$
\int_A  \frp_a(b)\, d\alpha(a)= \beta(b) 
,$$ 
or, precisely, for any subset
 $N\subset B$
$$
\int_A  \frp_a(N)\, d\alpha(a)= \beta(N)
. 
$$

Informally, we can consider 
 $\frp$ as a map  $A\to B$ spreading a point
$a\in A$  to a measure  $\frp_a$ on $B$ (or
spreading a point    $a$ 
along $B$).

\sm

{\sc Example.} Consider a bistochastic kernel  $\frp\in\Mar(A,B)$,
which is equal
 $\alpha\times \beta$. Then for all  $b$ we have
 $\frp_a(b)=\beta(b)$. 
The corresponding map uniformly 'spreads' any point
 $a$ along  $B$. For any morphism 
$\frq\in\Mar(B,C)$, we have
$$
\frq\circ (\alpha\times \beta)=\alpha\times\gamma
.
$$
For $\fro\in \Mar(Z,A)$ we have
$$
\qquad\qquad\qquad\qquad\qquad\qquad
(\alpha\times\beta)\circ \fro=(\zeta\times\alpha)
.
\qquad\qquad\qquad\qquad\qquad\qquad\square
$$

The product of $\frp\in\Mar(A,B)$, $\frq\in\Mar(B,C)$
can be regarded as an iterated spreading.
Precisely, let   $\frp_a$, $\frq_b$ be the corresponding systems
of conditional measures.
Then the measures $\frr_a(c)$ corresponding to $\frr=\frq\circ\frp$
are given by the formula 
$$
\frr_a(c)=\int_B  \frq_b(c)\,d\frp_a(b).
$$


{\bf \punct Markov operators.%
\label{ss:markov-operators}}
For a bistochastic kernel
 $\frp\in\Mar(A,B)$
 we define an operator
 $T(\frp)$ by the formula
$$
T(\frp) f(a)=\int_B f(b)\,d\frp_a(b)
.$$

\begin{proposition}
{\rm a)}
For each
 $p\in [1,\infty]$ the operator   $T(\frp)$
is bounded as an operator 
 $L^p(B)\to L^p(A)$.
 Moreover, its norm
 $\le 1$ for all $p$.
It also is continuous as an operator
 $L^{\infty-}(B)\to L^{\infty-}(A)$.

\sm

{\rm b)} The map $\frp\to T(\frp)$ is continuous as an operator
$L^{\infty_-}(B)\to L^{\infty_-}(A)$.

\sm

{\rm c)} The map $\frp\to T(\frp)$ is continuous
with respect to the weak operator topology%
\footnote{For definitions of weak and strong 
operator topologies, see 
 \cite{RS1}, \S VI.1}.
 
 \sm
 
{\rm d)} For $\frp\in\Mar(A,B)$, $\frq\in\Mar(B,C)$
we have
 $$
T(\frq\circ\frp)=T(\frp)T(\frq).
$$
\end{proposition}

{\sc Proof.} a)
Let $f\in L^p(A)$, $g\in L^q (B)$, where $\frac 1p+\frac 1q=1$.
Then 
\begin{multline*}
\Bigl|\int_A
T(\frp) g(a)\cdot f(a)\,d\alpha(a)\Bigr|=
\Bigl|
\int_{A\times B} f(a) g(b)\,d\frp(a,b)
\Bigr|
\le\\
\le \Bigl(\int_{A\times B} |f(a)|^p\,d \frp(a,b)  \Bigr)^{1/p}
\cdot
\Bigl(\int_{A\times B} |g(b)|^q\,d\frp(a,b)\Bigr)^{1/q}
\end{multline*}
(we applied the H\"older inequality, see \cite{Loe}, \S 9.3).
By the definition of the bistochastic kernel,
this equals
$$
\Bigl(\int_A |f(a)|^p\,d \alpha(a)  \Bigr)^{1/p}
\cdot
\Bigl(\int_B |g(b)|^q\,d\beta(b)\Bigr)^{1/q}=
\|f\|_{L^p}\cdot \|g\|_{L^q}.
$$
By the duality, we get the desired statement.
For $L^1$ a separate proof is necessary.
For a positive function
 $g$, the function  $T(\frp) g$ also is positive, moreover
$$
\int_B T(\frp) g(a)\,d\alpha(a)=
\int_{A \times B} g(a)d\frp(a,b)=\int_B g(b)\,d\beta(b) 
.$$
This implies the desired statement.

For $L^\infty$ the statement is obvious.

\sm

b) Our operator send the unit ball in $L^\infty(B)$ to a unit ball
 in
$L^\infty(A)$, moreover it is continuous
in the topology of  $L^1$. It remains to refer to
Lemma \ref{l:infty=p}.

c) Since the operators are uniformly bounded,
it suffices to verify the weak convergence on indicator functions%
\footnote{Let $M\subset A$ be a measurable subset,
we set $I_M(a)=1$ if $a\in M$, and $I_M(a)=0$ otherwise.%
\label{foo}}.
But
$$
\qquad\qquad\qquad
\int_A T(\frp) I_N(a)\, I_M(a)\,d\alpha(a)=\frp(M\times N). 
$$

{\rm d)} The statement is obvious for
absolutely continuous
kernels, the general case follows from separate continuity.
\hfill$\square$

\sm


\begin{figure}
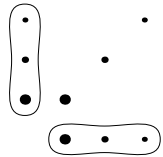


$$
\epsfbox{kvadratiki.2} \qquad\qquad \epsfbox{kvadratiki.6}
$$

The segment $A=[0,1]$, its partition $\sfX$ into   3 
pieces
and the quotient-space   $A/\sfX$.

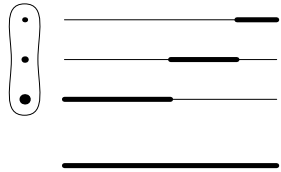
$$
\epsfbox{kvadratiki.1}
$$

Morphism
 $\frl[A;\sfX]$. On the figure, the product   $A/\sfX\times A$ 
 is a union of 3 horizontal segments. The measure 
 $\frl[A;\sfX]$ is a uniform measure
 on the union of free fat horizontal segments. 
 
 $$
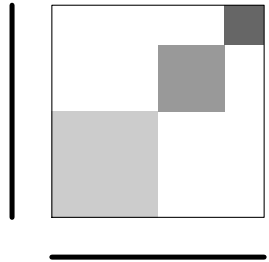
\epsfbox{kvadratiki.3}
$$

Morphisms  $\frt[A;\sfX]$. We have a uniform measure
on each subsquare    $\subset [0,1]\times [0,1]$.

$$
\epsfbox{kvadratiki.4}
$$

The identity morphism $A/\sfX\to A/\sfX$.

\caption{reference to Subsection \ref{ss:mar-conditional}. Morphisms
related to a partition.}

\end{figure}

\begin{figure}

$$
\epsfbox{kvadratiki.5} \qquad \qquad\qquad \epsfbox{kvadratiki.7}
$$

\bigskip

\qquad \qquad\qquad
$
\epsfbox{kvadratiki.8}
$

\caption{reference to Subsection \ref{ss:mar-conditional}.
The morphism $\frt[A/\sfX]\circ\frp\circ\frt[A/\sfX] $
is obtained from    $\frp$ 
by uniform spreading of the measure
 $\frp$ along each rectangle.
\newline
The morphism $\frm[A/\sfX]\circ\frp\circ\frl[A/\sfX] $ is obtained from   $\frp$ 
by concentration of measures
 $\frp$ on each rectangle.
}

\end{figure}

{\bf \punct Conditional expectations.%
\label{ss:mar-conditional}}
Let $\sfX: A=\cup_{r\in R} X_r$ be a measurable partition of  $A$, 
let
 $(R,\rho)$ be the quotient-space,
  $\pi: A\to R$
  the corresponding projection.
Consider a map 
$\xi:A\to A\times (A/\sfX)$ defined by the formula 
 $a\mapsto (a, \pi(a))$. Denote by 
$$
\frm[A;\sfX]\in \Mar(A,A/\sfX)
$$
the image of the measure  $\alpha$ under the map   $\xi$. Denote
$$
\frl[A;\sfX]:= \frm[A;\sfX]^\bigstar:\,\in\Mor(A/\sfX,A)
.$$

Define also the morphism
$$
\frt[A;\sfX]= \frl[A;\sfX]\circ \frm[A;\sfX]: A\to A
.
$$

Let us define such measures more explicitly.
The measure $\frm[A;\sfX]$ on $A\times A/\sfX$ is determined
from the equality
$$
\int_{A\times R} F(a,r) \,d\frm[A;\sfX](a,r)=
\int_{R}\biggl(\int_{X_r} F(a,r)\,d \xi_r(a)\biggr)\,d\rho(r)
.$$ 

The measure  $\frt[A;\sfX]$
is determined by
$$
\int_{A\times A} F(a_1,a_2) \, d\frt[A;\sfX]=
\int_R\biggl(\int_{X_r\times X_r } F(a_1,a_2) \,d\xi_r(a_1)\,d\xi_r(a_2)
\biggr)\,d\rho(r)
.$$

In notation of Subsection
\ref{ss:conditional-expectation},
$$
I[A;\sfX]=T\bigl(\frt[A;\sfX]\bigr),
\qquad
J[A;\sfX]=T\bigl(\frl[A;\sfX]\bigr),
\qquad
K[A;\sfX]=T\bigl(\frm[A;\sfX]\bigr)
.
$$

Now, let $\frp\in\Mar(A,B)$, let $\sfX: A=\cup X_i$, $\sfY: B=\cup Y_j$
be countable measurable partitions.
First, consider the measure
$$ 
\fru:=
\frm[B;\sfY]  \circ\frp\circ \frl[A;\sfX] \in \Mar(A/\sfX,B/\sfY)
.
$$
Both spaces
 $A/\sfX$, $B/\sfY$ are discrete. Therefore, the measure 
$\fru$ is determined by a matrix with non-negative elements,
this matrix equals
\begin{equation}
\fru_{ij}=\frp(X_i\times Y_j)
\label{eq:matrix-produce}
.
\end{equation}

Next, consider
$$ 
\frv:=
\frt[B;\sfY]  \circ\frp\circ \frt[A;\sfX] \in \Mar(A,B)
.
$$
This measure equals
$$
\frv(M\times N)
=
\sum_{i,j}\frac{\alpha(M\cap X_i)}{\alpha(M)}
\frac{\beta(N\cap Y_j)}{\beta(N)}  \frp(X_i\times Y_j)
,$$
where $M\subset A$, $N\subset B$ are measurable sets
of non-zero measure. 

\sm


{\bf\punct Definition of the product via approximations.%
\label{ss:mar-limit-def}}
Let $\sfX^{(1)}$, $\sfX^{(2)}$, \dots be a sequence of measurable
partitions of the set  $A$. We say 
that it is  {\it approximating},
if for any
 $p$ the partition 
 $\sfX^{(p+1)}$ is refinement of   $\sfX^{(p)}$
 and the sigma-algebra generated by all partitions coincides with
 the sigma-algebra of all measurable subsets of
 $A$.

Let $A$, $B$, $C$  be spaces with probability measures
and
 $\sfX^{(p)}$, $\sfY^{(q)}$, $\sfZ^{(r)}$ be approximating sequences
 of partitions of the spaces
 $A$, $B$, $C$ respectively. 

\begin{proposition}
\label{pr:mar-approximation}
The product $\frp\in\Mar(A,B)$, $\frq\in\Mar(B,C)$ equals
\begin{align*}
&\frq\circ \frp=\\
&\lim_{ i\, j,\, k\to\infty}
\frt[C;\sfZ^{k}]
    \circ\frq\circ \frt[B;\sfY^{l}] 
\circ
\frp\circ \frt[A;\sfX^{i}] =
\\
&\lim_{ i\, j,\, k\to\infty}
\frl[C;\sfZ^{k}]
\circ
\Bigl(  \frm[C;\sfZ^{k}]  \circ\frq\circ \frl[B;\sfY^{l}] \Bigr)
\circ
\Bigl( \frm[B;\sfY^{j}]  \circ\frp\circ \frl[A;\sfX^{i}] \Bigr)
\circ
\frm[A;\sfX^{i}]
\end{align*}
\end{proposition}

The products inside the brackets are elements of
$$
\Mar(A/\sfX^{(i)}, B/\sfY^{(l)})
\quad \text{and}\quad \Mar( B/\sfY^{(l)}, C/\sfZ^{(k)}).
$$
They are 
matrices of form
 (\ref{eq:matrix-produce}).
 Product of two brackets 
 (an element of $\Mar(A/\sfX^{(i)},C/\sfZ^{(k)})$)
 is calculated as a product of matrices
 as 
 (\ref{eq:psi-delta-phi}). Adding 
 outermost factors
 we get a measure
 on $A\times C$ and after this pass to the limit.

\sm

{\sc Proof.}
Let us pass to the corresponding Markov operators. The sequences 
$T\bigl( \frm[C;\sfZ^{k}] \bigr)$,  
$T\bigl(\frl[B;\sfY^{l}] \bigr)$, $T\bigl(\frl[A;\sfX^{i}] \bigr)$
 {\it strongly} converge to $1$.
 Let three strongly convergent sequences of operators
 are given, 
$A_i\to A$, $B_j\to B$, $C_k\to C$. Then
$$
A_iB_jC_k f-ABC f= A_iB_j (C_k-C)f+ A_i(B_j-B)C f+(A_i-A)BC f,
$$
and the right-hand side tends to 0 by norm as
 $(i,j,k)\to (\infty,\infty,\infty)$.
 This implies the desired statement.
\hfill $\square$


\section{Semiring of measures on  $\R^\times$}

\COUNTERS

This section is a preparation to the definition
of polymorphisms.


\sm

{\bf\punct Semiring $\cM^\tri$.%
\label{ss:tri}}
Denote by  $\cM^\tri$ the set of all positive measures
 $\mu$ on $\R^\times$
 satisfying the conditions
$$
\int_{\R^\times} d\mu(t)<\infty,\qquad
\int_{\R^\times}t\, d\mu(t)<\infty
.$$

Obviously, $\mu$, $\nu\in\cM^\tri$ implies $\mu+\nu\in \cM^\tri$.
We also introduce on  $\cM^\tri$
 the usual convolution $(\mu,\nu)\mapsto \mu*\nu$:
$$
\int_{\R^\times} f(t)\,d\mu*\nu(t)=\int\int_{\R^\times\times \R^\times} f(s_1s_2)\,d\mu(s_1)\,d\nu(s_2)
.$$

Obviously, $\cM^\tri$ is closed with respect to the convolution. 
Indeed,
\begin{equation}
\int\limits_{\R^\times} t^u d\mu*\nu(t)= 
\iint\limits_{(\R^\times)\times (\R^\times)} s_1^us_2^u\,d\mu(s_1)\,d\nu(s_2)=
\int\limits_{\R^\times} s_1^ud\mu(s_1) \cdot \int\limits_{\R^\times} s_2^ud\nu(s_2)
.
\label{eq:fff}
\end{equation}
Substituting
 $u=0$ and $u=1$ we get $\mu*\nu\in\cM^\tri$.

Next. we define an involution
$\mu\mapsto\mu^\bigstar$ in $\cM^\tri$ by the formula 
$$
\mu^\bigstar(t)=t^{-1}\mu(t^{-1}) 
.$$
In other words,
\begin{equation}
\int_{\R^\times} f(t)\,d\mu^\bigstar(t)=\int_{\R^\times} t f(t^{-1})\,d\mu(t)
.
\label{eq:mtri-conditions}
\end{equation}
If $\mu\in\cM^\tri$, then
$\mu^\bigstar\in\cM^\tri$; also
$(\mu*\nu)^\bigstar =\mu^\bigstar*\nu^\bigstar$.

We say that a sequence
 $\mu_j\in\cM^\tri$ converges to  $\mu\in \cM^\tri$,
 if for any bounded function
 $f(t)$ on $\R^\times$ the following convergences hold
$$
\int f(t)\,d\mu_j(t)\to \int f(t)\,d\mu(t),
\qquad 
\int t f(t)\,d\mu_j(t)\to \int t f(t)\,d\mu(t).
$$
In other words, we require a weak convergence
(see, e.g.,  \cite{Loe}, Sect. 12.1)
of two sequences of measures
$\mu_j\to \mu$, $t\mu_j\to t\mu$.


\sm

{\bf\punct Mellin transform.%
\label{ss:m-mellin}}
{\it Mellin transform} of a measure 
 $\mu\in\cM^\tri$  is
\begin{equation}
\Phi_\mu(u):=
\int_{\R^\times} t^u \,d\mu(t),\qquad\text{where $u=v+i w\in\C$.}
\label{eq:Phi}
\end{equation}

{\sc Remark.}
Pass to a variable
 $s:=\ln t$. Then $\nu(s)=\mu(\ln t)$ is a measure on 
 $\R$,  conditions  (\ref{eq:mtri-conditions})
take the form
$$
\int_\R  d\nu(s)<\infty, \qquad \int_\R e^s \nu(s)<\infty
.
$$
The function
 $\Phi(u)$ is the characteristic function    (Fourier transform)
of the measure $\nu$. This topic is quite standard 
  (see, e.g.,  \cite{Lin}), however I do not see 
  an appropriate reference for lemmas given below. 
\hfill$\square$

\begin{proposition}
\label{pr:M-mellin}
{\rm a)} For any $\mu\in\cM^\tri$ the function
$\Phi_\mu$ is uniformly continuous in the strip 
\begin{equation}
\Pi:\,
0\le v\le 1 \qquad -\infty< w<\infty
\end{equation}
and holomorphic in the open strip
 $0<\Re u< 1$.
 
\sm

{\rm b)} The functions  $\Phi_\mu(u)$ are positive definite,
i.e., for any  $u_1$, \dots, $u_n$,
satisfying  $0\le \Re u_j\le 1/2$ and any $z_1$, \dots, $z_n\in \C$
\begin{equation}
\sum_{l,m\le n}
\Phi(u_l+\ov u_m)z_l\ov z_m\ge 0
.
\label{eq:positve-definite} 
\end{equation}

{\rm c)} The functions $\Phi_\mu$ satisfy the following estimate
$$
|\Phi_\mu(v+iw)|\le 
\Phi(0)^{1-v} \Phi(1)^{v} 
.
$$
In particular,
 $\Phi_\mu(u)$ is bounded in the strip   $0\le \Re u\le 1$.
\end{proposition}

{\sc Proof.}
 a) The convergence of integral  (\ref{eq:Phi})
 is obvious.
Let us prove uniform continuity:
$$
|\Phi_\mu(u)-\Phi_{\mu}(u')|\le
\int_{\R^\times} |t^u-t^{u'}|\,d\mu(t)
$$
We split the integral as a sum of  integrals over segments 
$t<1/A$, $1/A\le t\le B$, $t>B$. First,
$$
\int_{t>B} |t^u-t^{u'}|\,d\mu(t)\le \int_{t>B} 2t\,d\mu(t)
.
$$
For sufficiently large 
 $B$ the integral is as small as desired.
Similarly, we estimate the integral over
 $t< 1/A$:
$$
\int_{t<1/A} |t^u-t^{u'}|\,d\mu(t)\le \int_{t<1/A} 2\,d\mu(t)
.
$$
Next, fix
  $A$, $B$, 
\begin{multline*}
\int_{1/A\le t\le B} |t^u-t^{u'}|\,d\mu(t)=
\int_{1/A\le t\le B} t^{\Re u}|t^{u'-u}-1|\,d\mu(t)
\le 
\\
\le
\int_{1/A\le t\le 0} |t^{u'-u}-1|\,d\mu(t)
+
\int_{0 < t\le B} t\,|t^{u'-u}-1|\,d\mu(t)
\end{multline*}
For small $|u'-u|$ a value $|t^{u'-u}-1|$ is small $[1/A,B]$.

A proof of
 b) is usual,
$$
0\le
\int_{\R^\times} \bigl|\sum_k z_k t^{u_k}\bigr|^2 \,d\mu(t)=
\sum_{k,l}\Phi(u_k+\ov u_l)\,z_k \ov z_l.
$$ 
 
  To prove   c)
 we apply H\"older inequality, 
$$\qquad
\Bigl| \int t^{u+iw}\, d\mu(t)\Bigr|
\le \Bigl(\int 1^{1/(1-v)} \,d\mu(t)\Bigr)^{1-v}
\cdot \Bigl(\int |t^{v+iw}|^{1/v} \,d\mu(t)\Bigr)^v
\qquad\qquad\square
$$

\begin{proposition}
\label{pr:M-mellin-inverse}
Let $\Phi(u)$ be a bounded positive definite
function continuous in the strip
 $0\le \Re u\le 1$ holomorphic in the open strip.
Then
$\Phi(u)$ is a Mellin transform of some measure
 $\mu\in \cM^\tri$. 
\end{proposition}

{\sc Proof.} By the Paley--Wiener theorem, see \cite{Her},  Theorem 7.4.2,
the function
$\Phi$ is a Fourier transform of a certain distribution
  $\nu(s)$ that is contained in the Schwartz space
  on $\R$. Applying the Bochner Theorem   
(see, e.g., \cite{Loe}, \S 15.1) to the function  $\Phi(iw)$,
we get that
   $\nu(s)$ is a finite positive measure.
Applying the Bochner Theorem to  $\Phi(1+iw)$
we get   that $e^s\cdot \nu(s)$
is a finite  measure. Passing to the variable  $t=e^s$, 
we get the desired statement. \hfill$\square$

\begin{proposition}
\label{pr:m-mellin-product}
{\rm a)}
$
\Phi_{\mu*\nu}(u)=\Phi_\mu(u)\Phi_\nu(u)
$.

{\rm b)}
$\Phi_{\mu^\bigstar}(u)=\Phi_\mu(1-u)$
\end{proposition}

{\sc Proof.} a) is evident, it was proved
by calculation (\ref{eq:fff}); b) 
also is obvious.
 \hfill $\square$


\sm

{\bf \punct Convergence of characteristic functions.%
\label{ss:convergent-characteristic}}

\begin{proposition}
\label{pr:m-mellin-convergence}
If $\mu_j$ converges to $\mu$ in $\cM^\tri$, then $\Phi_{\mu_j}(u)$
converges to $\Phi_\mu(u)$ uniformly in each
rectangle 
$0\le v\le 1$, $-A\le w\le B$.
\end{proposition}

The pointwise convergence is obvious, a proof of
uniform convergence coincides with a standard proof,
see
\cite{Loe}, 13.2.C.
\hfill $\square$

\begin{proposition}
\label{pr:m-mellin-continuity}
{\rm a)} Let $\mu_j$, $\mu\in \cM^\tri$. If
\begin{equation}
\Phi_{\mu_j}(iw)\to \Phi_{\mu_j}(iw),\qquad
\Phi_{\mu_j}(1+iw)\to \Phi_{\mu_j}(1+iw)
\label{eq:pointwise}
\end{equation}
pointwise, then
 $\mu_j$ converges to $\mu$.

\sm

{\rm b)} Let $\mu_j\in \cM^\tri$. Let the sequence 
$
\Phi_{\mu_j}(iw)$ converge point-wise to a function  $\Psi(iw)$
and
$\Phi_{\mu_j}(1+iw)$ converge pointwise to a function  $\Theta(1+iw)$.
If $\Psi(iw)$, $\Theta(1+iw)$ are continuous at the point
 $w=0$, then $\mu_j$
 converges to a certain
 $\mu\in\cM^\tri$, and
$\Phi(iw)=\Psi(iw)$,
$\Phi(1+iw)=\Theta(1+iw)$.
\end{proposition}

{\sc Proof.} Let us prove b). By the continuity theorem 
(see, e.g.,  \cite{Loe}, Theorem 15.2),
the sequence
 $\mu_j$ weakly converges to a certain
  measure  $\mu$ and $t\cdot \mu_j$
  weakly converges to a certain measure 
 $\nu$. Let $f(t)$   be a  continuous function with compact
 support.
  Then
\begin{multline*}
\int_{\R^\times} f(t)\,d\nu(t)=\lim_{j\to\infty} \int_{\R^\times} f(t) t\,d\mu_j(t)
=\\=
\lim_{j\to\infty} \int_{\R^\times} \bigl(t\,f(t)\bigr) \,d\mu_j(t)=
\int_{\R^\times}  \bigl(t\,f(t)\bigr)\,d\mu(t)
\end{multline*}
Hence,
 $\nu(t)=t\mu(t)$ and $\mu_j$ converges to  $\mu$ in the sense of $\cM^\tri$.
\hfil$\square$


\sm


\section{Polymorphisms. Basis definitions}

\COUNTERS

{\bf\punct Definition.%
\label{ss:polymorphisms-def}} Let $(A,\alpha)$, $(B,\beta)$
be Lebesgue measure spaces.
A {\it polymorphism} $A\zigzag B$  is a measure  $\frP$ on $A\times B\times \R^\times$
such that

\sm

$1^\circ$. The image of  $\frP$ under the projection   $A\times B\times \R^\times \to A$
coincides with $\alpha$;

\sm

$2^\circ$. The image of $t\cdot\frP$ under the projection 
 $A\times B\times \R^\times \to B$
coincides with $\beta$.

\sm

We denote the set of all polymorphisms by
 $A\zigzag B$ by $\Pol(A,B)$.


There is a well-defined associative product
$$
\Pol(A,B)\times\Pol(B,C)\to \Pol(A,C).
$$
A formal definition is given in
 Subs.\ref{ss:pol-product-2}. 
 Before this in Subsections \ref{ss:mar-pol}--\ref{ss:abs-cont-ker} we consider several simple special cases.

\sm


{\bf\punct A special case: the category   $\Mar$.%
\label{ss:mar-pol}}
Let $A$, $B$ be spaces with probability measures.
An element $\frp\in \Mar(A,B)$  can be regarded
as an element from 
 $\Pol(A,B)$, we simply consider an image
 of the measure  
$\frp$ under the embedding 
$$A\times B\to  A\times B\times \R^\times$$
defined by the formula
 $(a,b)\mapsto (a,b,1)$.

\sm


{\bf\punct A special case: $\cM^\tri$.%
\label{ss:tri-pol}}
Consider single-point spaces
 $A$, $B$, denote by  $\alpha$, $\beta$ their measures. 
Then a polymorphism 
$A\zigzag B$ is a measure on $\R^\times$ satisfying
$$
\int_{\R^\times} d\frP(t)=\alpha,
\qquad 
\int_{\R^\times} t d\frP(t)=\beta
.
$$
A product of polymorphisms 
 $\frP:A\zigzag B$, $\frQ:B\zigzag C$ 
 coincides with the convolution of measures
 in $\cM^\tri$:
 $$
 \frQ\circ\frP=\frac 1\beta \frQ*\frP.
 $$

\sm


{\bf\punct Special case: discrete spaces.%
\label{ss:pol-discrete}} Consider discrete spaces
 $A$, $B$, let $a_i$, $b_j$ be their points,  
$\alpha_i$, $\beta_j$ the measures of these points.
A measure  $\frP$ on $A\times B\times \R^\times$
can be regarded as a matrix, whose matrix elements
are non-negative measures on
$\frp_{ij}\in\cM^\tri$,
these measures must satisfy the condition
\begin{align}
\sum_i \int_{\R^\times}t\,d\frp_{ij}=\beta_j;
\label{eq:discrete-identities-1}
\\
\sum_j \int_{\R^\times}d\frp_{ij}=\alpha_i
\label{eq:discrete-identities-2}
.
\end{align}

For $\frP\in\Pol(A,B)$, $\frQ\in\Pol(B,C)$, 
their product 
  $\frR=\frQ\circ\frP$ is defined by
 \begin{equation}
 \frr_{ik}=\sum_j \frac 1{\beta_j} \frq_{jk}*\frp_{ij}
, 
\label{eq:discrete-product}
\end{equation}
where $*$ denotes the convolution in   $\cM^\tri$.
Actually, we have matrices, whose elements are measures
   $\in\cM^\tri$,
see (\ref{eq:psi-delta-phi}).


{\bf\punct   A special case: absolutely continuous
kernels.%
\label{ss:abs-cont-ker}} 
Let $p:A\times B\to \cM^\tri$ be a measurable function. 
We define a measure
 $\frP$ on $A\times B\times \R^\times$
 in the following way.
For measurable subsets $M\subset A$, $N\subset B$, $K\subset \R$ we assume
$$
\frP(M\times N\times K):=\int_M\int_N p(a,b)(K)\,d\beta(b)\,d\alpha(a)
.$$
If
\begin{align}
\int_B\int_{\R^\times} d p(a,b)(t)\,d\beta(b)=
1, \qquad\text{a.s.}, 
\label{eq:abs1}
\\
\int_A \int_{\R^\times}t dp(a,b)(t)\,d\alpha(a)=1
,\qquad\text{a.s.}, 
\label{eq:abs2}
\end{align}
then $\frP$ is a polymorphism. In this case we say that   $\frP$ 
is
{\it  absolutely continuous}.

Obviously, a polymorphism
 $\frP$ is absolutely continuous if the projection
 $\frP$ to $A\times B$ is a measure
 absolutely continuous with respect to
  $\alpha\times \beta$.

\sm

{\sc Remark.} This class of polymorphisms 
includes the objects of the previous subsection. 
For a matrix   $\frp_{ij}$ the function $p$ is given by
$$
\qquad\qquad\qquad\qquad\qquad\qquad
p(a_i\times b_j)=\frac{\frp_{ij}}{\alpha_i\beta_j}
.
\qquad\qquad\qquad\qquad\qquad\qquad\square
$$

Let $\frP:A\zigzag B$, $\frQ:B\zigzag C$ be absolutely continuous
polymorphisms,
$p$, $q$ the corresponding    $\cM^\tri$-valued functions.
We define a function  $r:A\times C\to\cM^\tri$ by the formula
$$
r(a,c)=\int_B p(a,b)*q(b,c)\,d\beta(b)
.$$

\begin{lemma}
{\rm a)} $r(a,c)\in \cM^\tri$ a.s.

\sm

{\rm b)} $r$  determines a polymorphism  $A\zigzag C$.
\end{lemma}

{\sc Proof.}
For a proof of
  a) we write an integral
\begin{multline*}
\int_B\int_A\int_{\R^\times} t\, d (p(a,b)*q(b,c))\,d\alpha(a)\,d\beta(b)
=\\
=\int_B \Bigl( \int_{\R^\times} t\,d q(b,c)(t))\Bigr)
\int_A\Bigl(\int_{\R^\times} t\, d p(a,b)(t)\Bigr) \,d\alpha(a)\,d\beta(b)
=\\=
\int_B \Bigl( \int_{\R^\times} t\,d q(b,c)(t))\Bigr)
\,d\beta(b)=1
\end{multline*}
and change an order of integration
to
$\int_A\int_B\int_{\R^\times}$
By the Fubini theorem
the integral
$$
\int_B\int_{\R^\times} t\, d (p(a,b)*q(b,c))\,d\beta(b)
$$
converges a.s.
Next, we repeat the same argument for the integral
$$
\int_B\int_C\int_{\R^\times} \, d (p(a,b)*q(b,c))\,d\gamma(c)\,d\beta(b)
.$$

b) is verified by a straightforward calculation. 
\hfill $\square$


\sm

{\bf\punct Definition of the product.%
\label{ss:pol-product-2}}
Now we  define the product of morphisms.
Below in 
Subsection \ref{ss:pol-product-3} we give another
(may be, more transparent) definition.
Also, Theorem
\ref{th:pol-mellin-product} can be used a definition.

Let $\frP\in \Pol(A,B)$. 
For any measurable subsets
 $M\subset A$, $N\subset B$ we have  a measure $\frp[M\times N]\in \cM^\tri$
defined as the image of 
 $\frP$ under the projection
$$
M\times N\times \R^\times\to \R^\times.
$$
In this space we can regard
 $\frP$ as a $\cM^\tri$-valued measure  $\frp[\cdot]$ on $A\times B$.

\begin{lemma} 
\label{l:for-product}
{\rm a)} Let $\frP\in\Pol(A,B)$.  For any measurable subset
  $M\subset A$ there is a system of measures 
  $\frp_{M,b}(t)$, 
where $b$ ranges in $B$,  on  $\R^\times $,
defined for almost all
  $b\in B$,
  such that for any measurable 
  $N\subset B$ we have 
\begin{equation}
\frp[M\times N]=\int_N \frp_{M,b}\,d\beta(b)
.
\label{eq:frp}
\end{equation}

{\rm b)} Let $\frQ\in\Pol(B,C)$. For any measurable subset 
 $K\subset C$ there is a system of measures 
$\frq_{b,K}(t)$ on $\R^\times$ such that
for any measurable subset
 $N\subset B$ we have
\begin{equation}
\frq[N\times K]
=\int_N \frq_{b,K}\, d\beta(b)
.
\label{eq:frq}
\end{equation}
\end{lemma}

{\sc Proof.} 
a) Denote by $\frS_M$ the restrictions of the measure
  $t\frP$ to  $M\times B\times \R^\times$. Consider images of
  $\frS_M$ under the projections
$$
M\times B\times \R^\times\stackrel{p}{\longrightarrow} B\times\R^\times \stackrel{q}{\longrightarrow} B
.
$$
The measure $q(p(t\frS_M))$ is dominated by  $\beta$. 
Therefore there are well-defined conditional measures 
$\sigma_{M,b}(t)$ on fibers of the projection $B\times \R^\times\to B$, such that
 $$
 p(\frS_M)(N)=\int_N \sigma_{M,b}(N)\,d\beta(b).
 $$
By the construction, $\sigma_{M,b}(\R^\times) \le 1$ (since
$\sigma_{A,b}(\R^\times)=1$).
We define measures $\frp_{M,b}$ as
$$
\frp_{M,b}:=t^{-1} \sigma_{M,b}(t).
$$

b) It suffices to consider  images of  $\frQ$ under the projections
$$
\qquad\qquad\qquad\qquad\qquad
B\times K\times \R^\times\to B\times \R^\times\to B
\qquad\qquad\qquad\qquad\qquad
.\square
$$

\sm

To each subset
 $M\times K\subset A\times C$ we assign the measure
\begin{equation}
\frr[M\times K]=
\int_B \frq_{b,K}*\frp_{M,b}\,d\beta(b)\,\,\in\cM^\tri
\label{eq:frr}
\end{equation}
and we get a  $\cM^\tri$-valued measure
on $A\times C$.

\begin{lemma}
\label{l:countable-additive}
{\rm a)}
$\frr$  is a sigma-additive    $\cM^\tri$-valued measure on  $A\times C$.

\sm

{\rm b)}  The measure $\frr$ determines a polymorphism 
$A\zigzag C$.
\end{lemma}

Lemma is proved in next subsection.

\begin{lemma}
For absolutely continuous kernels this product coincides with
the product defined above.
\end{lemma}

 {\sc Proof.} Let $p$, $q$, $r$ be the same as in
Subsection \ref{ss:abs-cont-ker}. Then
 $$
 \frp_{M,b}=\int_M p(a,b)\,d\alpha(a),
 \qquad
 \frq_{b,K}=\int_K q(b,c)\,d\gamma(c).
 $$
Therefore
\begin{multline*}
\frr(M\times K)=\int_B\int_{M\times K}  p(a,b)*q(b,c)\,d\alpha(a)\,d\gamma(c)\,d\beta(b)
=\\=
\int_{M\times K} r(a,c)\,d\alpha(a)\,d\gamma(c).
\end{multline*}

\begin{theorem}
\label{th:assoc-product}
The product
 $\Pol(A,B)\times \Pol(B,C)\to\Pol(A,C)$ is associative,
 i.e. for any measure spaces 
 $A$, $B$, $C$, $D$ and any
$\frP\in\Pol(A,B)$, $\frQ\in\Pol(B,C)$, $\frT\in\Pol(C,D)$
we have
$$
(\frT  \circ \frQ) \circ \frP = \frT  \circ (\frQ \circ \frP)
.
$$
\end{theorem}

Proof is given below in Subsection
\ref{ss:associativity-proof}.

\sm


{\bf\punct  Proof of lemma \ref{l:countable-additive}.}
First, we need a more  precise information about functions
 $\frp_{M,b}$ and
$\frq_{b,K}$ defined in Lemma
 \ref{l:for-product}.

\begin{lemma}
\label{l:lfc1}
{\rm a)} $\frp_{M,b}\in \cM^\tri$ a.s. on $b\in B$.
\begin{equation}
\text {\rm b)}\qquad \qquad \qquad \qquad 
 \int_B\int_{\R^\times} d\frp_{M,b}(t)\,d\beta(b)= \alpha(M).\qquad \qquad \qquad \qquad 
\label{eq:l1}
\end{equation}
\begin{equation}
\text {\rm c)} \qquad \qquad \qquad 
\int_{\R^\times} t \,d\frp_{M,b}(t)\le 1 \quad\text{\rm  for almost all $b\in B$.}
 \qquad \qquad \qquad 
\label{eq:l2}
\end{equation}
и
\begin{equation}
\int_{\R^\times} t d\frp_{A,b}(t)= 1 
\label{eq:l201}
\end{equation}

{\rm d)} If   $\alpha(M_j)$ tends to  0, then
\begin{equation}
\int_{B}\int_{\R^\times} t\cdot d\frp_{M_j,b}(t)\,d\beta(b)\,\to 0
.
\label{eq:l202}
\end{equation}
\end{lemma}

{\sc Proof.} Statements  b), c) follow from the same arguments
as Lemma 
\ref{l:for-product}.
By (\ref{eq:l1}) the measures $\frp_{M,b}$ are finite for almost
all  $b$.
 By (\ref{eq:l2}) they are contained in 
 $\cM^\tri$.

The projection of the measure 
 $t\cdot\frP$ to $A$ is a probabilistic measure 
 absolutely continuous with respect to
 $\alpha$. The statement  d) is rephrasing of this fact.
\hfill $\square$

\sm
 
 Next, we formulate a similar lemma for measures
 $\frq_{b,K}$.

\begin{lemma}
\label{l:lfc2}
{\rm a)} $\frq_{b,K}\in \cM^\tri$ for almost all  $b\in B$.
\begin{equation}
\text{\rm b)} \qquad \qquad \qquad \qquad 
\int_B\int_{\R^\times} t\,d\frq_{b,K}(t)\, d\beta(b)= \gamma(K)
\qquad \qquad \qquad 
\label{eq:l3}
\end{equation}
\begin{equation}
\text{\rm c)} \qquad \qquad 
\int_{\R^\times} \,d\frq_{b,K}(t)\le 1 \quad\text{\rm for almost all $b\in B$.}
\qquad \qquad \qquad 
\label{eq:l4}
\end{equation}
и
\begin{equation}
\int_{\R^\times} \,d\frq_{b,C}(t)= 1 \quad\text{\rm for almost all $b\in B$.}
\qquad \qquad \qquad 
\label{eq:l401}
\end{equation}

{\rm d)} If $\gamma(K_j)\to 0$, then
$$
\int_B
\int_{\R^\times} \,d\frq_{b,K_j}(t)\,d\beta(b)\to 0
$$
\end{lemma}

A proof is the same.

\sm

{\sc Proof of Lemma  \ref{l:countable-additive}.a.}
If $M_1$, $M_2$ are disjoint, then  
$$\frp_{M_1,b}+\frp_{M_2,b}=\frp_{M_1\cup M_2,b}.$$
 By 
(\ref{eq:frr}) this implies finite additivity.

For a proof of sigma-additivity, we take a decreasing chain of subsets
 $M_1\supset M_2\supset \dots$ in $A$,
such that $\alpha(M_j)\to 0$:
\begin{multline*}
\int_{\R^\times} d\frr[M_j\times K](t)=\int_B 
\Bigl(\int_{\R^\times} d\frp_{M_j,b}(t)\Bigr)
\cdot
\Bigl(\int_{\R^\times} d\frq_{b,K}(t)\Bigr)\,d\beta(b)
=\\ =
\alpha(M_j)\int_B \int_{\R^\times} d\frq_{b,K}(t)\,d\beta(b)\to 0
,
\end{multline*}
here we applied
 (\ref{eq:l1}), (\ref{eq:l4}).
Further,
\begin{multline*}
\int_{\R^\times} t\cdot d\frr[M_j\times K](t)=
\int_B \Bigl(\int_{\R^\times} t\cdot d\frp_{M_j,b}(t)\Bigr)
\cdot
\Bigl(\int_{\R^\times}t\cdot d\frq_{b,K}(t)\Bigr)\,d\beta(b)
=\\=
\gamma(K)
\int_B \int_{\R^\times} t\cdot d\frp_{M_j,b}(t)
\,d\beta(b)\to 0
,\end{multline*}
here we applied
 (\ref{eq:l3}), (\ref{eq:l202}).

 Thus, $\frr[M_j\times K]\to 0$ in $\cM^\tri$.
 
 Let now $L_j$ be mutually disjoint subsets in 
  $A$, $M_j=\cup_{i\le j} L_i$. In virtue of proved above,
  we get the following limit pass 
 in $\cM^\tri$:
 $$
 \frr[\cup_{i=1}^\infty L_i \times K]=\lim_{j\to\infty} \frr[\cup_{i=1}^j L_i \times K]=
 \sum_{j=1}^\infty \frr[L_j\times K].
 $$

 Similarly, we consider a decreasing chain 
 $K_1\supset K_2\supset \dots$,
such that  $\gamma(K_j)\to 0$.

In the same way, we prove that
 $\frr[M\times K_j]\to 0$ as
$\gamma(K_j)\to\infty$.
\hfill $\square$

\sm

{\sc Proof of Lemma  \ref{l:countable-additive}.b.}
\begin{multline*}
\int_{\R^\times}d\frr[M\times C](t)
=
\int_B 
\Bigl(\int_{\R^\times} d\frp_{M,b}(t)\Bigr)
\cdot
\Bigl(\int_{\R^\times} d\frq_{b,C}(t)\Bigr)\,d\beta(b)
=\\=
\int_B 
\int_{\R^\times} d\frp_{M,b}(t)
\,d\beta(b)=\alpha(M)
,\end{multline*}
here we applied
 (\ref{eq:l401}) and (\ref{eq:l1}). Next,
\begin{multline*}
\int_{\R^\times}t\cdot d\frr[A\times K](t)
=
\int_B 
\Bigl(\int_{\R^\times}t\cdot  d\frp_{A,b}(t)\Bigr)
\cdot
\Bigl(\int_{\R^\times}t\cdot  d\frq_{b,K}(t)\Bigr)\,d\beta(b)
=\\=
\int_B 
\int_{\R^\times}t\cdot  d\frq_{b,K}(t)\,d\beta(b)=
\gamma(K)
,
\end{multline*}
here we applied (\ref{eq:l201}) и (\ref{eq:l3}).
\hfill$\square$


\sm

{\bf \punct Involution.%
\label{ss:pol-involution}} Let $\frP\in \Pol(A,B)$.
We define  $\frP^\bigstar\in \Pol(B,A)$ being the measure 
$t^{-1}\frP(a,b,t^{-1})$ regarded as a measure 
$B\times A\times \R^\times$. 

\begin{lemma}
\label{l:pol-involution}
For $\frP\in\Pol(A,B)$, $\frQ\in\Pol (B,C)$, we have 
$$
(\frQ\circ \frP)^\bigstar=\frP^\bigstar\circ \frQ^\bigstar
.
$$
\end{lemma}

{\sc Proof.} Multiplying $\frP^\bigstar\circ \frQ^\bigstar$
we get in 
 (\ref{eq:frr}) the expression 
$$
\int_B (t^{-1} \frq_{b,K})*(t^{-1}\frp_{M,b})\,d\beta(b)=t^{-1}\frr[M\times K](t^{-1})
.$$


\sm

{\bf\punct Convergence.%
\label{ss:pol-convergence}}
Let $\frP^{(j)}$, $\frP\in\Pol(A,B)$. We say that a sequence
  $\frP_j$ converges to 
 $\frP$ if for any measurable subsets 
$M\subset A$, $N\subset B$ we have the convergence
$$
\frp^{(j)}(M\times N) \to \frp(M\times N)
$$ 
in the sense of
 $\cM^\tri$.

\begin{theorem}
\label{th:pol-separately}
The $\circ$-product is separately continuous.
\end{theorem}

{\sc Proof.} We keep the notation of  
Subs.\ref{ss:pol-product-2}.
Since we have an involution, it is sufficient to prove 
one-side continuity.
 Let $\frP^{(j)}\to \frP$.
This means that for any measurable $N\subset B$
 functions
$\frp_{M,b}$  satisfy the condition
\begin{equation}
\int_N \frp^{(j)}_{M,b} \,d\beta(b) \quad\text{converges to}\quad 
\int_N \frp_{M,b} \,d\beta(b)\quad\text{in $\cM^\tri$}.
\label{eq:mmm}
\end{equation}
Also
$$
\int_B \frp^{(j)}_{M,b} \,d\beta(b) \le \frp[A\times B].
$$
We intend to show that
$$
\frr^{(j)}[M\times K]\to \frr[M\times K] \quad\text{in $\cM^\tri$}.
$$
It is sufficient to verify point-wise convergence
of Mellin transforms on the lines
 $u=iw$, $u=1+iw$. We have
\begin{multline}
\int_{\R^\times} t^{iw}d\bigl(\frr^{(j)}[M\times K]-\frr[M\times K]\bigr)\,d\beta(b)
=\\=
\int_B
\Bigl[\int_{\R^\times} t^{iw} d\frp^{(j)}_{M,b}(t)- \int_{\R^\times} t^{iw} d\frp_{M,b}(t)\Bigr]
\cdot\Bigl\{ \int_{\R^\times} t^{iw} d\frq_{b,K}(t)\Bigr\}\,d\beta(b)
.
\label{eq:convergence-long}
\end{multline}
The factor   $\{F(b)\}$ in curly brackets
is a bounded function,
see (\ref{eq:l4}). The factor 
$[G^{(j)}(b)-G(b)]$  in square brackets 
is contained in  $L^1(B)$ and its    $L^1$-norm is bounded by a constant
  $2\alpha(M)$ (by (\ref{eq:l1})).
On the other hand, for any  $N$
 $$
 \int_N [G^{(j)}(b)-G(b)]\,d\beta(b)\to 0\quad\text{as $j\to\infty$.}
 $$ 
For any  $\epsilon>0$ we take a function   $F_\epsilon(b)$
taking only finite number of values such that
$$
\ess | F(b)-F_\epsilon(b)|<\epsilon
.
$$
Then
$$
\int_B [G^{(j)}(b)-G(b)]\cdot F_\epsilon(b)\,d\beta(b)\to 0\qquad
\text{for $j\to\infty$}.
$$
On the other hand
$$
\int_B [G^{(j)}(b)-G(b)]\cdot (F_\epsilon(b)-F(b))\,d\beta(b)\le \epsilon \cdot 2\alpha(M).
$$
Therefore
 (\ref{eq:convergence-long}) tends to $0$. 

Further, 
\begin{multline}
\int_{\R^\times} t^{1+iw}d\bigl(\frr^{(j)}[M\times K]-\frr[M\times K]\bigr)\,d\beta(b)
=\\=
\int_B
\Bigl[\int_{\R^\times} t^{1+iw} d\frp^{(j)}_{M,b}(t)- \int_{\R^\times} t^{1+iw} d\frp_{M,b}(t)\Bigr]
\cdot\Bigl\{ \int_{\R^\times} t^{1+iw} d\frq_{b,K}(t)\Bigr\}\,d\beta(b)
.
\label{eq:convergence-long-2}
\end{multline}

Again, denote the expression in square brackets by
$[G_j(b)-G(b)]$,
the expression in curly brackets by $F(b)$.
Now
 $F\in L^1(B)$ by (\ref{eq:l3}),
and
 $ [G_j(b)-G(b)]\le 2$ by (\ref{eq:l2}),
 moreover, for any
 $N\subset B$ the following convergence holds
$$
\int_N[\dots]\,d\beta(b)\to 0\quad \text{as $j\to\infty$,}
$$
i.e., we have a convergence in the sense of the space
$L^{\infty_-}(B)$.
Consider a function  $F_{\epsilon,\delta}$ taking only
finite number of values 
such that
$$
\ess | F(b)-F_{\epsilon,\delta}(b)|<\epsilon
\qquad\text{on a set of measure $>\beta(B)-\delta$.}
$$
It is sufficiently obvious that
$$
\int_B [G_j(b)-G(b)] \cdot F_{\epsilon,\delta}(b)\,d\beta 
\to 0\qquad \text {for $j\to\infty$,}
$$
and
$$
\int_B [G_j(b)-G(b)] \cdot (F(b)- F_{\epsilon,\delta}(b))\,d\beta 
$$
is small for small
 $\epsilon$ и $\delta$.
Therefore the expression
(\ref{eq:convergence-long-2}) tends to  0.
\hfill$\square$

\sm



{\bf\punct  Proof of Theorem  \ref{th:assoc-product}.
Associativity of the product.%
\label{ss:associativity-proof}}
The set of absolutely continuous polymorphisms
 $A\zigzag B$ is dense  $\Pol(A,B)$.
Obviously, the product of absolutely continuous 
polymorphisms is associative.
On the other hand, a product of polymorphisms is separately continuous.
 \hfill$\square$


\sm

{\bf \punct Definition of the product in the terms 
of discrete approximations.%
\label{ss:pol-product-3}}
Let us return to the definition of Subsection
\ref{ss:mar-conditional}. 
For a countable partition
  $\sfX$ of the set
$A$ we  define  morphisms
 $$
 \frl[A;\sfX]:A/\sfX\zigzag A, \qquad\frm[A;\sfX]:A\zigzag A/\sfX,\qquad
\frt[A;\sfX]:A\to A
$$
 as above   (recall that
$\Mar(A,B)\subset\Pol(A,B)$).
Consider countable
partitions  $\sfX:A=\cup X_i$, $\sfY:B=\cup Y_j$.

\begin{lemma}
{\rm a)} For $\frP: A\zigzag B$ a morphism
$$
\frm[A;\sfY]
\circ
\frP
\circ
\frl[A;\sfX]:\, A/\sfX\zigzag B/\sfY
$$
is given by
$\cM^\tri$-valued matrix  $\frp_{ij}=\frp[X_i\times Y_j]$.

\sm

{\rm b)} The measure
$$
\frt[A;\sfY]
\circ
\frP
\circ
\frt[A;\sfX]:\, A\zigzag B
$$
is determined by the following rule:
its restriction to
 $X_i\times Y_j\times\R^\times\subset A\times B\times \R^\times$
coincides with 
$$\frac 1{\alpha(X_i)\beta(Y_j)}
\cdot
\alpha\times\beta\times \frp_{ij}
.$$
\end{lemma}

A verification is straightforward.
In any case, the statement follows from
Theorem
 \ref{th:pol-mellin-product} proved below.

\sm

For measure spaces 
 $A$, $B$, $C$ consider approximating sequences
 (see Subs. \ref{ss:mar-limit-def})
 of countable partitions
$\sfX^{(i)}$, $\sfY^{(j)}$, $\sfZ^{(k)}$.

\begin{proposition}
\label{pr:pol-approximative}
The product of polymorphisms 
 $\frP:A\zigzag B$, $\frQ:B\zigzag C$ 
 is given by the formula
\begin{multline}
\frQ\circ \frP=
\lim_{ i,\, j,\, k\to\infty}
\frt[C;\sfZ^{(k)}]
    \circ\frQ\circ \frt[B;\sfY^{(l)}] 
\circ
\frP\circ \frt[A;\sfX^{(i)}]
 =
\\
=
\lim_{ i,\, j,\, k\to\infty}
\frl[C;\sfZ^{(k)}]
\circ
\Bigl(  \frm[C;\sfZ^{(k)}]  \circ\frQ\circ \frl[B;\sfY^{(l)}] \Bigr)
\circ\\ \circ
\Bigl( \frm[B;\sfY^{(l)}]  \circ\frP\circ \frl[A;\sfX^{(i)}] \Bigr)
\circ
\frm[A;\sfX^{(i)}]
.
\label{eq:triple-limit}
\end{multline}
\end{proposition}

The expressions in big brackets are polymorphisms 
of countable sets, their product is evaluated in the way
described above
 (\ref{eq:discrete-product}).

\sm

A proof of the proposition is given
in Subsection
\ref{ss:proof-approximation}.

\sm

{\sc Proof.} Notice, that a reference to the separate continuity
allows to claim that
 $\frQ\circ \frP$ coincides with the iterated limit   
$$
\frQ\circ \frP=
\lim_{ i\to\infty} \lim_{ j\to\infty} \lim_{ k\to\infty}\bigl(\dots).
$$
But in
 (\ref{eq:triple-limit}) we have triple limit.
\hfill $\square$

\sm


{\bf\punct The group $\Gms(A)$.%
\label{ss:pol-gms}} Let $A$  be a space with continuous measure.
For $g\in\Gms(A)$ consider the map 
$\frI_g :A\to A\times A\times \R^\times$ defined by
$$
a\mapsto \bigl(a,g(a),g'(a)\bigr)
.
$$
Denote by
 $\frI[g]$ the image of measure 
$\alpha$ under this map.

\begin{figure}

$$
\epsfbox{krivulina.1}
$$
The map $y=x+\frac 1n \sin nx$ of the segment  $[0,2\pi]$ to itself.
$$
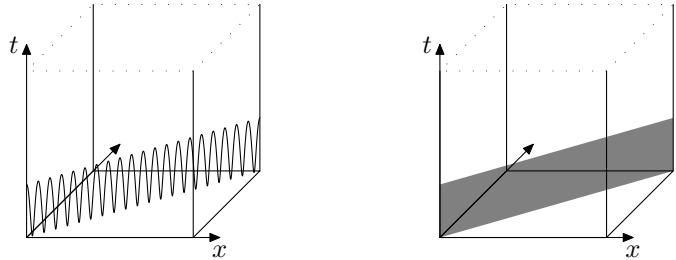
\epsfbox{krivulina.2}
\qquad\qquad\qquad
\epsfbox{krivulina.3}
$$

The image of the segment $[0,1]$ in
$[0,2\pi]\times[0,2\pi]\times \R^\times$ is
 an oblate helical line. The limit as
 $n\to\infty$ is a (non-uniform) measure supported by
 the rectangle 
 $x=y$, $0<t\le 2$. 

\caption{reference to Subsection \ref{ss:pol-gms}}

\end{figure}

\begin{proposition}
\label{pr:gms-pol}
{\rm a)} $\frI[g]\in \Pol(A,A)$

\sm

{\rm b)} The map $g\mapsto \frI(g)$ is a homomorphism.
\end{proposition}

This is obvious.

\begin{theorem}
\label{th:gms-dense}
Let
 $A$ be a space with continuous measure. 
Then the group  
$\Gms(A)$ is dense in  $\Pol(A,A)$.
\end{theorem}

{\sc Proof.} Fix $\frP\in \Pol(A,A) $. Consider a finite partition
$\sfX:A=\cup X_j$ of the space $A$. Denote
$\frp_{ij}=\frp(X_i\times X_j)$.
Consider a subdivision
 $X_i=\cup Y_{ij}$ of each  $A_i$ such that
$$
\alpha(Y_{ij})=\int_{\R^\times}  \,d\frp_{ij}(t)$$
Consider another subdivision
 $X_i=\cup Z_{ij}$ 
 such that
$$
\alpha(Z_{ij})=\int_{\R^\times} t \,d\frp_{ij}(t)
.$$
For any pair
 $(i,j)$ consider the map 
$Y_{ij}\to Z_{ij}$, whose Radon-Nikodym derivative
is distributed as
 $\frp_{ij}$. 
 Uniting the maps
 $X_{ij}\to Y_{ij}$ we get an element 
$g[\sfX]$ of the group $\Gms(A)$.

Next,  we consider an approximating sequence of partitions
 $\sfX^{(p)}$ and get the sequence 
  $g[\sfX^{(p)}]\in\Gms(A)$, which converges to  $\frP$.
\hfill $\square$

\sm



\section{Mellin--Markov transform}

\COUNTERS


\sm

 
{\bf\punct Mellin--Markov transform.%
\label{ss:markov-mellin}} Let $u=v+iw$  range in the strip 

\begin{equation}
\Pi:\,
0\le v\le 1 \qquad -\infty< w<\infty
.
\end{equation}
Denote $p=1/(1-v)$, $q=1/v$.
Let
$\frP\in\Pol(A,B)$ be a polymorphism.
Consider the following bilinear form on
 $L^{1/(1-v)}(A)\times L^{1/v}(B)$
\begin{equation}
S_{u}(\frP;f,g)=S_{v+iw}(\frP;f,g)= \iiint_{A\times B\times \R^\times} f(a)g(b)t^{v+iw}\,d\frP(a,b,t)
.
\label{eq:bilinear}
\end{equation}

\begin{lemma}
\label{l:holder}
 $$
|S_{v+iw}(\frP;f,g)|\le \|f\|_{1/(1-v)} \cdot \|g\|_{1/v}
.
$$
\end{lemma}

The lemma, in particular, implies the continuity of the bilinear form
on the product of Banach spaces for all
 $u\in \Pi$. We need  an improvement of the statement for  
for the spaces
 $L^{\infty_-}$. 

\begin{lemma}
\label{l:infty-}
Let a sequence of bounded functions
 $f_j$ converges to  $f$ in $L^{\infty_-}(A)$,
 and $g\in L^1(B)$. Then
$S_{1+iw}(\frP;f_j,g)$  converges to $S_{1+iw}(\frP;f,g)$.
\end{lemma}

{\sc Proof of Lemma \ref{l:holder}.} We apply the H\"older
inequality 
 (see \cite{Loe}, \S 9.3) and the definition of polymorphisms
\begin{multline*}
|S_{v+iw}(f,g)|\le
\biggl(
\iiint_{A\times B\times \R^\times} |f(a)|^{1/(1-v)}\,d\frP(a,b,t) 
\biggr)^{1-v}
\times\\\times
\biggl(
\iiint_{A\times B\times \R^\times} |g(b)t^{v+iw}|^{1/v}\,d\frP(a,b,t)
\biggr)^{v}
=\\=
\biggl(\int_A |f(a)|^{1/(1-v)}\,d\alpha(a)\biggr)^{1-v}
\cdot 
\biggl(
\int_{ B} |g(b)|^{1/v}\,d\beta(b)
\biggr)^{v}
.
\qquad \square
\end{multline*}

{\sc Proof of Lemma \ref{l:infty-}.}
Without loss of generality, we can set
 $f=0$. 
Denote by  $L$ the set, where
$|f_j(a)|\ge \epsilon $. Split integral on two summands, via
$(A\setminus L)\times B\times \R^\times$ and $L\times B\times \R^\times$,
 \begin{multline*}
\Bigl| \iiint\limits_{(A\setminus L)\times B\times \R^\times} f_j(a)g(b)t^{1+iw}\,d\frP(a,b,t)\Bigr|
\le
\epsilon  \iiint\limits_{(A\setminus L)\times B\times \R^\times} |g(b)|\,t\,d\frP(a,b,t)
\le \\\le \epsilon  \iiint_{A\times B\times \R^\times} |g(b)|\,t\,d\frP(a,b,t)
=\epsilon  \int_B |g(b)|\,d\beta(b).
 \end{multline*}
Estimate another summand
$$
\Bigl| \iiint\limits_{L\times B\times \R^\times} f_j(a)g(b)t^{1+iw}\,d\frP(a,b,t)\Bigr|
\le \ess\limits_{a\in A} |f_j(a)|
 \iiint\limits_{L\times B\times \R^\times} |g(b)|\,t\,d\frP(a,b,t)
$$
Notice that
\begin{equation}
 \iiint\limits_{L\times B\times \R^\times} t\,d\frP(a,b,t)
 \label{eq:mal}
\end{equation}
is small if
 $\alpha(L)$ is small. Indeed,
 denote by 
 $\sigma$ the projection of the measure  $t\,d\frP(a,b,t)$ to $A$.
Then $\sigma(A)=\beta(B)$, the measure $\sigma$
is absolutely continuous with respect to  $\alpha$.
Indeed, take in $A$ a subset   $S$ of zero measure.
Then the measure of $S\times B\times \R^\times$ with respect to $\frP$ is $0$.
Therefore $\frP(S\times B\times (0,y])=0$,  hence
$(t\cdot\frP)(S\times B\times (0,y])=0$,
 so $(t\cdot \frP)(S\times B\times \R^\times)=0$.

Consequently (by the absolute continuity 
of the Lebesgue integral), $\alpha(L)\to 0$
implies $\frP(L\times B\times \R^\times)\to 0$.
Applying the absolute continuity of the Lebesgue integral again,
we get that (\ref{eq:mal})
tends to 0 as
 $\alpha(L)\to 0$.
\hfill $\square$

\sm

As an immediate corollary of these lemmas we get the following 
theorem

\begin{theorem}
\label{th:pol-mellin-def}
{\rm a)}
For any  $u=v+iw\in \Pi$ there exists a linear operator 
$$
T_u(\frP): L^p(B)\to L^p(A), \qquad \text{where $p=\frac 1{v}$}
,$$
satisfying
\begin{equation}
\int_A \bigl(T_u(\frP)g\bigr)(a)\, f(a)\,d\alpha(a)=S_u(\frP;f,g)
\label{eq:T-def}
\end{equation}
for all $f\in L^q(A)$.

\begin{equation}
{\rm b)}\qquad\qquad\qquad\qquad\qquad
\|T_u(\frP)\|_{L^p}\le 1.
\qquad\qquad\qquad\qquad\qquad
\end{equation}

\end{theorem}

We call the map
 $u\mapsto T_u(\frP)$ by the  {\it Mellin--Markov transform}
 of the polymorphisms  $\frP$.
 
 \sm

{\sc Remark.} Let  $s<1/v<r$. Then  $L^r\subset L^{1/v}\subset L^s$,
i.e., $T_u$ is bounded as an operator
 $L^r\to L^s$. In particular all the operators $T_u$
 are bounded as operators $L^{\infty_-}\to L^1$.
\hfill $\square$
 
 \sm
 
{\sc Remark.} The same form determines a
dual operator 
$T_u(\frP)': L^{1/(1-v)}(A)\to L^{1/(1-v)}(B)$.
\hfill $\square$

\begin{lemma}
 $S_u(\frP^\bigstar;f,g)=S_{1-u}(\frP;g,f)$.
\end{lemma}

{\sc Proof.} We substitute
 $t\mapsto t^{-1}$
to (\ref{eq:bilinear}).
\hfill $\square$

\sm

As a corollary we get:

\begin{proposition}
The operator $T_u(\frP^\bigstar)$ is dual to 
$T_{1-u}(\frP)$.
\end{proposition}
 

{\bf\punct Direct definition of the Markov--Mellin transform.%
\label{ss:pol-product-1}}
First, we reformulate the definition of polymorphisms.
Fix a polymorphism
 $\frP:A\zigzag B$. Consider a map 
$A\times B\times \R\to A$.
For  $a\in A$ consider the conditional (probabilistic) measure
  $\frP_a(b,t)$ om $B\times\R^\times$.
  Next, consider the map
$B\times \R^\times \to B$.
Denote the image of the measure  $\frP_a(b,t)$ by
 $\frP_a(b)$. By  $\frP_{a,b}(t)$ we denote the conditional
 measures on fibers. In other words,
\begin{equation}
\iiint\limits_{A\times B\times \R^\times} F(a,b,t)\,d\frP(a,b,t)
=\int\limits_A\biggl(\int\limits_B\biggl(\int\limits_{\R^\times}
 F(a,b,t)\,d\frP_{a,b}(t)\biggr)\,d\frP_a(b)  \biggr)\,d\alpha(a)
 \label{eq:fubbinization}
 .
\end{equation}
Now we can define a polymorphism in the terms  of two systems
of conditional measures
$\frP_a(b)$, $\frP_{a,b}(t)$. These measures
are probabilistic and satisfy to the integral identity
corresponding to the condition 
$2^\circ$ for polymorphisms (see Subs. \ref{ss:polymorphisms-def}):
\begin{equation}
\int\limits_A\biggl(\int\limits_B\biggl(\int\limits_{\R^\times}
 t\, g(b)\,d\frP_{a,b}(t)\biggr)\,d\frP_a(b)  \biggr)\,d\alpha(a)
 =\int_B g(b)\,d\beta(b).
  \label{eq:integral-identity-0}
\end{equation}
This holds for
 $g\in L^1(B)$.  The identity can be written also as
\begin{equation}
\int\limits_A\biggl[ \biggl(\int\limits_{\R^\times}
 t\,d\frP_{a,b}(t)\biggr)\cdot\frP_a(b)  \biggr]\,d\alpha(a)=
 \beta(b)
 \label{eq:integral-identity}
 .
\end{equation}
In the square brackets there is a product of an integrable function and a measure.

\begin{theorem}
\label{th:mellin-def-2}
For $\frP:A\zigzag B$ and $g\in L^1(B)$ 
the following equality holds
\begin{equation}
T_u(\frP)g(a)=\int_B\int_{\R^\times}
t^u g(b)\,d\frP_{a,b}(t)\,d\frP_a(b)
.
\label{eq:escho}
\end{equation}
\end{theorem}

{\sc Proof.}
For an operator (\ref{eq:escho}),
$$
\int_A T_u(\frP)g(a) \,f(a)\,d\alpha(a)
=
\int_A
\int_B\int_{\R^\times}
t^u f(a)\, g(b)\,d\frP_{a,b}(t)\,d\frP_a(b)
\,d\alpha(a)
.
$$
By the definition of our conditional measures
 (see (\ref{eq:fubbinization})),
$$d\frP_{a,b}(t)\,d\frP_a(b)\,d\alpha(a)=d\frP(a,b,t),$$
an we get $S_u(f,g)$.
\hfill $\square$

For absolutely continuous kernels the formula is more transparent.
Let $p:A\times B\to\cM^\tri$ 
be the same function as in
\ref{ss:abs-cont-ker}.

\begin{proposition}
$$
T_u(\frP) g(a)=
\int_B\int_{\R^\times} t^u g(b) \,dp(a,b)(t)\,d\beta(b).
$$
\end{proposition}

This and the following statements are obvious.

\begin{proposition}
{\rm a)} For polymorphisms  $\frP\in\Mar(A,B)$ the operators $T_u(\frP)$  
coincide with Markov operators defined 
in Subs.
{\rm\ref{ss:markov-operators}}.

\sm

{\rm b)} For $g\in\Gms(A)$ these operators coincide with operators
 $T_u(g)$ defined by the formula
{\rm(\ref{eq:t1p})}.

\sm

{\rm c)} For single-point spaces 
$A$, $B$
the function  $u\mapsto T_u$ 
coincides with characteristic function
  $u\mapsto\Phi(u)$ discussed in 
{\rm\S 4.}
\end{proposition}

Notice one's more corollary

\begin{proposition}
 The operators $T_{iw}$ are continuous as operators
 $L^{\infty_-}(B)$ to
$L^{\infty_-}(A)$.
\end{proposition}

{\sc Proof.} We intend to apply
 Lemma \ref{l:infty=p}. Let
$\ess |g|=1$ and $\|g\|_{L^1}$ be small:
Estimate the absolute value of $T_{iw}(\frP)g$
$$
|T_{iw}(\frP) g(a)|
=\Bigl|\int_B\int_{\R^\times}
t^{iw} g(b)\,d\frP_{a,b}(t)\,d\frP_a(b)\Bigr|
\le \int_B\int_{\R^\times}
|g(b)|\,d\frP_{a,b}(t)\,d\frP_a(b)
.
$$
Integrating the last expression over 
 $A$ we get
$$
\int_A[\dots]\,d\alpha(a)=
\int_{A\times B\times\R^\times} |g(b)|\,d\frP(a,b,t).
$$
The projection of the measure
$\frP(a,b,t)$ to $B$  is absolutely continuous
with respect  $\beta$ and finite
(see proof of Lemma \ref{l:infty-}),
i.e. has the form $h(b)\,d\beta(b)$  with integrable  $h$.
We come to
\begin{equation}
\int_B |g(b)|\, h(b)\, d\beta(b).
\label{eq:lll}
\end{equation}
Let $S\subset B$  be the set of points $b$, where $|g(b)|>\epsilon$.
We split the integral into two summands,
$$
\int_{B\setminus S} |g(b)|\, h(b)\, d\beta(b)\le
 \epsilon \int_{B\setminus S} \, h(b)\, d\beta(b)
 \le 
 \epsilon \int_{B} \, h(b)\, d\beta(b)=\epsilon \alpha(A), 
$$
$$
\int_{S} |g(b)|\, h(b)\, d\beta(b)\le \ess |g(b)|
\cdot \int_{S}  h(b)\, d\beta(b)
.
$$
By the absolute continuity of the Lebesgue integral, we get that the last expression
is small for small
 $\beta(S)$. \hfill $\square$


{\bf\punct Holomorphy of matrix elements.}
The function
 $u\mapsto T_u(\frP)$ is holomorphic in the following
 sense.

\begin{lemma}
\label{l:strip}
{\rm a)}
For fixed
 $f\in L^\infty(A)$, $g\in L^\infty(B)$ the function
$u\mapsto S_u(\frP;f,g)$ is continuous in the strip
 $\Pi$ and holomorphic in the open strip function.

 \sm
 
 {\rm b)} If $f$, $g\in L^\infty$ are non-negative, then
 the function  $u\mapsto S_u(\frP;f,g)$ is positive definite in the strip
 $\Pi$.
\end{lemma}

{\sc Proof.} a) The integral 
$$
\frac \partial {\partial u}
S_u(\frP;f,g)=\iiint_{A\times B\times \R^\times} f(a)g(b)t^u\ln(t)\,d\frP(a,b,t)
$$
converges for
 $0<v<1$ and is dominated by a convergent 
 integral in a neighborhood of any
 point.

\sm

b) The measure $f(a)g(b)\,d\frP(a,b,t)$ is contained in $\cM^\tri$.
\hfill $\square$

\sm

To be complete, we present a more precise statement about
matrix elements.

\begin{lemma}
\label{l:matrix-elements}
 For $f\in L^{r}(A)$, $g\in L^{s}(B)$ the function   $S_u(\frP;f,g)$ 
 is continuous in the strip
\begin{equation}
\frac 1r\le \Re u\le 1-\frac 1s
\label{eq:strip-2}
\end{equation}
and holomorphic in the corresponding open strip.
\end{lemma}

{\sc Proof.} The bilinear form $S_{v+iw}(\frP;f,g)$
is continuous on $L^{1/(1-v)}(A)\times L^{1/v}(B)$. Therefore it is
continuous on
$L^s(A)\times L^r(B)$ for $s>1/(1-v)$, $r>1/v$.
\hfill $\square$

\sm


{\bf\punct Characterization of the image and inversion.}

\begin{theorem}
\label{th:pol-inversion}
{\rm a)} 
Let $u\mapsto T_u$ be a function in the strip
  $\Pi$, 
  taking values in the space of bounded operators
 $L^{\infty_-}(A)\to L^1(B)$,
 such that

\sm

{\rm i)}  For  positive
 $f\in L^\infty(A)$, $g\in L^\infty(B)$
 matrix elements
$$
u\mapsto \phi_{f,g}(u)= \int_A T_u g(a)\, f(a)\,d\alpha(a)
$$
are continuous and bounded
for
 $u\in\Pi$ and holomorphic in the open strip;

\sm

{\rm ii)} For non-negative   $f$, $g$  functions  $\phi_{f,g}(u)$
are positive definite in
$\Pi$;

\sm

{\rm iii)} 
$\phi_{1,1}(0)=1$, $\phi_{1,1}(1)=1$.

\sm

Then there exists a unique polymorphism
 $\frP\in\Pol(A,B)$ such that  $T_u=T_u(\frP)$. 

\sm

{\rm b)} The polymorphism  $\frP$ is determined by the condition:
$$
\int_{\R^\times} t^u \,d\frp[M\times N](a,b,t)=\phi_{I_M, I_N}(u)
$$
for any measurable $M\subset A$, $N\subset B$, {\rm(}where $I_M$ 
denotes an indicator function of a set,
see footnote {\rm\ref{foo})}.
\end{theorem}

{\sc Proof.}
The function $\phi_{I_M, I_N}(u)$ is positive definite
and bounded in the strip   $\Pi$, 
therefore it is a characteristic function of a measure
$\frp:=\frp[M\times N]$. 

Obviously, for disjoint sets
 $M_1$, $M_2$ we have
$$
\frp\bigl[(M_1\cup M_2)\times N\bigr]=\frp[M_1\times N]+\frp[M_2\times N]
.$$
Further, let measurable sets
 $M_1$, $M_2$, $\dots\subset A$ be mutually disjoint.
Then $I_{\cup M_j}=\sum I_{M_j}$ in topology of   $L^{\infty_-}$.
Therefore the sequence
 $\phi_{I_{M_1\cup\dots\cup M_j},I_N}(u)$ converges pointwise to
$\phi_{I_{\cup_j M_j},I_N}(u)$.
By  Proposition \ref{pr:m-mellin-continuity}  we get
$$
\sum_j \frp[M_j\times N]=\frp\bigl[( \cup M_j)\times N\bigr]
.$$
The same arguments prove the similar equality for $M\times N_j$.
This implies that the
 $\cM^\tri$-valued measure on $A\times B$ is sigma-additive.

In the virtue of condition
  iii),  this measure is a polymorphism.
\hfill $\square$

\sm

{\bf \punct Convergence.}

\begin{theorem}
 \label{th:pol-mellin-convergence}
{\rm a)} If $\frP_j$ converges to $\frP$, then for any  $u\in \Pi$, $\Re u>0$,
the operators  $T_u(\frP_j):L^{1/v}(B)\to L^{1/v}(A)$ weakly converge
 to $T_u(\frP)$. In the case  $v=0$ we have a weak convergence
 in the space of operators  
 $L^{\infty_-}(B)\to L^{\infty_-}(A)$.
 
\sm

{\rm b)}
Let $\frP_j$, $\frP\in\Pol(A,B)$. Let
$$T_u(\frP_j):L^{\infty_-}(B)\to L^1(A)$$
weakly converges to
$T_u(\frP)$ for any   $u\in \Pi$. Then $\frP_j$  converges to $\frP$.
It is sufficient to require weak convergence on the lines
 $u=iw$ and $u=1+iw$.
\end{theorem}

{\sc Proof.}
a) It is sufficient to prove that
\begin{equation}
\int_A T_u(\frP_j) I_N(a)\,I_M(a)\,d\alpha(a)
\quad\text{converges to} \int_A T_u(\frP) I_N(a)\,I_M(a)\,d\alpha(a)
\label{eq:PP}
\end{equation}
for any $M\subset A$, $N\subset B$. We write this as a convergence
\begin{equation}
\int_{\R^\times} t^u\, d\frp_j[M\times N](t)\to \int_{\R^\times} t^u\, d\frp[M\times N](t)
\label{eq:pp}
\end{equation}
of Mellin transforms in
 $\cM^\tri$. Now we can refer to Proposition 
 \ref{pr:m-mellin-convergence}.

\sm

b) In virtue of Proposition \ref{pr:m-mellin-continuity} 
we have convergence (\ref{eq:pp}), this is equivalent to (\ref{eq:PP}).
\hfill $\square$


\sm

{\bf\punct Product.}

\begin{theorem}
\label{th:pol-mellin-product}
{\rm a)} For any $\frP\in\Pol(A,B)$, $\frQ\in\Pol(B,C)$,
the following identity holds
\begin{equation}
 T_u(\frP)T_u(\frQ)=T_u(\frQ\circ \frP)
.
\label{eq:mellin-product}
\end{equation}
\end{theorem}

 {\sc Proof.} The statement is obvious for absolutely continuous 
 kernels. 
 In the virtue of separate continuity of the product of polymorphisms
  and the  weak separate continuity of products of operators,
the statement holds for arbitrary polymorphisms.

\sm


{\bf\punct Proof of Proposition  \ref{pr:pol-approximative}.%
\label{ss:proof-approximation}} 
It is easy to see that
$$
T_u(\frt[A,\sfX])= I[A;\sfX]
.
$$
If a sequence 
 $\sfX^{(i)}$ is approximating, then    $I[A;\sfX^{(i)}]$
 strongly converges to 1.
A product of strongly convergent sequences
of operators strongly converges
(see the proof of Proposition  \ref{pr:mar-approximation}). 
Therefore the sequence
\begin{multline*}
T_u\Bigl(
\frt[C;\sfZ^{(k)}]
    \circ\frQ\circ \frt[B;\sfY^{(l)}] 
\circ
\frP\circ \frt[A;\sfX^{(i)}]
\Bigr)=\\=
I[C;\sfZ^{(k)}]\,\,
    T_u(\frQ)\,\, I[B;\sfY^{(l)}] \,\,
T_u(\frP)\,\, I[A;\sfX^{(i)}]
\end{multline*}
converges to
$$
T_u(\frQ)T_u(\frP)=T_u(\frQ\circ \frP)
$$
for $i$, $j$,  $k\to\infty$. Therefore
the triple sequence  
$\frt[C;\sfZ^{(k)}]
    \circ\frQ\circ \frt[B;\sfY^{(l)}] 
\circ
\frP\circ \frt[A;\sfX^{(i)}]$
strongly converges to
 $\frQ\circ \frP$.
\hfill $\square$



\tt

 Math.Dept., University of Vienna
 
 ITEP (Moscow)
 
 Moscow State University, MechMath
 
 URL:www.mat.univie.ac.at/$\sim$neretin

\end{document}